\documentclass[12pt]{article}

\usepackage{amsmath, amssymb, amsthm, mathrsfs,physics}
\usepackage{hyperref}  
\usepackage{enumitem}  
\usepackage{geometry}  
\numberwithin{equation}{section}

\newtheorem{theorem}{Theorem}[section]
\newtheorem{lemma}[theorem]{Lemma}
\newtheorem{proposition}[theorem]{Proposition}
\newtheorem{corollary}[theorem]{Corollary}

\theoremstyle{definition}
\newtheorem{definition}[theorem]{Definition}

\theoremstyle{remark}
\newtheorem{remark}[theorem]{Remark}
\newtheorem{example}[theorem]{Example}

\allowdisplaybreaks

\title{Discrete frames of non-uniform shifts in frequency domain}

\author{ Hari Krishan Malhotra, Manisha Chhillar, and Lalit Kumar Vashisht\\
\small{}
}
\date{June 18, 2025}

\begin{document}
	
	\maketitle
	
	\begin{abstract}
Non-uniform frames play a significant role in the stable analysis and decomposition of signals (functions) when the indexing set may not be a group (or, uniform lattice).  We analyze frame conditions for a  collection of matrix-valued  functions obtained by non-uniform shifts. We give necessary and sufficient conditions for the existence of matrix-valued  discrete Bessel sequence over non-uniform displacement parameters in terms of the Fourier transform of window functions. We also present perturbation results for matrix-valued non-uniform discrete frames. 
	\end{abstract}
\baselineskip15pt	
	\section{Introduction}
Based on the pioneering work by Gabor  \cite{Gabor},  Duffin and Schaeffer introduced the concept of frames in the published paper  \cite{DS}. They tried to find completeness of a family of complex exponential in the Lebesgue space  $L^2(-a, a)$, where $a >0$. Young in his book \cite{Young} reviewed frames and their variants. Let $H$ be a  separable (finite or infinite-dimensional) Hilbert space with respect to an inner product $\langle\cdot, \cdot\rangle$. Let $\mathbb{J}$ be a countable index set. A sequence   $\{x_n\}_{n\in\mathbb{J}} \subset H$ is  a \emph{discrete frame}  for $H$ if the following inequality holds for some positive real numbers $A$ and $B$:
\begin{align}\label{defn1}
A\|x\|^2\leq  \sum_{n\in\mathbb{J}} | \langle x, x_n \rangle |^2 \leq
B\|x\|^2,  \  x \in H.
\end{align}
Numbers $A$ and $B$  are called the \emph{lower frame bound} and \emph{upper frame bounds} of the frame $\{x_n\}_{n\in\mathbb{J}}$, respectively. The frame bounds of a frame are not unique. $\{x_n\}_{n\in\mathbb{J}}$ is a Bessel sequence with \emph{Bessel bound} $B$ if only upper inequality holds in \eqref{defn1}. If $\{x_n\}_{n\in\mathbb{J}}$ is a Bessel sequence, then the map $V: \ell^2(\mathbb{J})\rightarrow H$ defined  by
$V(\{\alpha_n\}_{n\in\mathbb{J}})=\sum_{n \in \mathbb{J}} \alpha_n x_n,  \ \{\alpha_n\}_{n\in\mathbb{J}} \in \ell^2(\mathbb{J})$,
is called the \emph{pre-frame operator} or the \emph{synthesis
operator} associated with $\{x_n\}_{n\in\mathbb{J}}$. The \emph{analysis operator} of $\{x_n\}_{n\in\mathbb{J}}$ is the Hilbert-adjoint operator $V^*$ of the pre-frame operator $V$. Note that  $V^*: H \rightarrow \ell^2(\mathbb{J})$ is given by
$V^{*}(x)=\{\langle x, x_n \rangle \}_{n\in\mathbb{J}}, \ x \in H$.     The \emph{frame operator} is the composition $S:VV^*: H \rightarrow H$ is given by $S(x)=\sum_{n\in\mathbb{J}}\langle x, x_n\rangle x_n$ $x \in H.$ If $\{x_n\}_{n\in\mathbb{J}}$ is a frame for $H$, then the frame operator $S$ is bounded, linear, positive, and invertible on $H$. Thus, each $x \in H$ can be decomposed as a series (not necessarily unique): $x = \sum_{n \in \mathbb{J}} \langle x, S^{-1}x_n\rangle x_n$.  In recent years, frames have emerged as a mathematical tool in the analysis of signals \cite{Heil1, Heil2, Kova07}, time-frequency analysis \cite{Groch}, algorithms related frames \cite{BHan2}, data science \cite{Zpelle}, atomic decompositions \cite{GroI}, dynamical sampling \cite{ABK}, approximation theory \cite{Larson}, distributed signal processing \cite{deep1, deep2, deep3}, etc. For application of frames of different structure, for example, frames of wavelets,  we refer to texts  by Krivoshein,  Protasov and   Skopina \cite{KPSA}, Novikov, Protasov and    Skopina \cite{NPS}, Han \cite{BHan2}, Heil \cite{Heil2} and many references therein.

The Fourier transform of a signal  in time domain gives information in the frequency domain. In \cite{zZhang}, Zhang studied support of scaling functions in terms of the Fourier transform. He characterized supports of Fourier transforms of scaling functions associated with  multiresolution analyzes (MRA). In \cite{Hari0}, two authors of this paper studied scaling functions of nonuniform multiresolution analysis in the frequency domain. There is huge literature on signal representation in frequency domain \cite{Hughes, KPSA, Tura}.  Frames are useful in stable analysis and decomposition of signals \cite{Groch, BHan2, Heil1, Heil2}.

In many applications where the data to be processed are samples of matrix-valued functions, for example, video images, multi-spectral images \cite{AZ}. In order to analyze signals in vector-valued signal spaces, Xia and Suter\cite{XS} introduced and studied vector-valued wavelets and vector-valued multiresolution analysis in the space $L^2(\mathbb{R}, \mathbb{C}^{n\times n})$. Xia and  Suter \cite{XS}  showed that the component functions in vector-valued wavelets generate multiwavelets. They also observed  in  \cite{XS} that multiwavelets can be constructed by taking certain linear combinations of the known scalar-valued wavelets. Xia in \cite{Xia1} studied  orthonormal matrix-valued wavelets and matrix Karhunen-Lo\`{e}ve expansion. Antolin and Zalik \cite{AZ} proposed matrix-valued wavelets which provide stable analysis and decompositions of matrix-valued signals. They characterized matrix-valued wavelet sets in terms of matrix-valued multiresolution analysis. Recently, frame conditions for matrix-valued  systems with different structure were studied in \cite{JindII, Jlk2,  Jlk4}.

On the other hand, in the study of  commuting self-adjoint partial differential operators, Fuglede \cite{CSPDG}  introduced the concept of spectral pairs. Using the concept of spectral pairs, Gabardo and Nashed  in \cite{Gaba1} introduced and studied  non-uniform multiresolution analysis in the signal space $L^2(\mathbb{R})$. Also see, \cite{Gaba2} for Cohen's condition for non-uniform multiresolution analysis. In non-uniform multiresolution analysis, the corresponding translation set may not be a group.
This concept is useful in  input/output models based on shifts of data points, where shifts may not be uniform in nature.
 Gabardo and Nashed \cite{Gaba1} gave the following definition:
\begin{definition} \cite{Gaba1} Let $\Omega \subset \mathbb{R}$ be measurable and $\Lambda \subset \mathbb{R}$ a countable subset. If the collection $\left\{|\Omega|^{-\frac{1}{2}} e^{2 \pi i \lambda \cdot} \chi_{\Omega}(\cdot)\right\}_{\lambda \in \Lambda}$ forms complete orthonormal system for $L^{2}(\Omega)$, where $\chi_{\Omega}$ is the  indicator function on $\Omega$ and $|\Omega|$ is Lebesgue measure of $\Omega$, then the pair $(\Omega, \Lambda)$ is a spectral pair.
\end{definition}
\begin{example} \cite{Gaba1} \label{exp12}Let $N \in \mathbb{N}, r$ be an odd integer coprime to $N$ such that $1 \leq r \leq 2 N-1$, and let $\Lambda=\left\{0, \frac{r}{N}\right\}+2 \mathbb{Z}$ and $\Omega=\left[0, \frac{1}{2}\left)\bigcup\left[\frac{N}{2}, \frac{N+1}{2}\Big)\right.\right.\right.$. Then, $(\Omega, \Lambda)$ is a spectral pair.
\end{example}
Frame properties of a non-uniform system obtained by applying the modulation operator and translation operator on a family of functions in a discrete signal space can be found in \cite{Hari02}. Frames of non-uniform wavelets were studied by two authors in \cite{Hari01}.  Motivated by above work, we study frame conditions for a family of non-uniform shift of the form:
\begin{align*}
   \left\{\mathcal{D}_{2 N \lambda} \mathbf{f}_{j}: \lambda \in \Lambda, j=1,2,\cdots,p \right\},
\end{align*}
 where $\mathcal{D}_{(.)}$ is a displacement operator  on $\ell^{2}\left(\Lambda,\mathcal{M}_{n}\right)$, and $\mathbf{f}_{j}$ are \emph{envelop sequences} in $\ell^{2}\left(\Lambda,\mathcal{M}_{n}\right)$ (see Definition \ref{defnM1} below).

\subsection{Overview of Main Results}
The main results are given in Section 2 and Section 3. In Sec.2, Theorem \ref{thmM1} gives a sufficient condition for the existence of discrete matrix-valued Bessel sequence, for the discrete matrix-valued signal space $\ell^{2}\left(\Lambda,\mathcal{M}_n\right)$,  in the  of nonuniform shifts in terms of Fourier transforms of signals. A necessary condition for discrete matrix-valued Bessel sequence of nonuniform shift in the frequency domain is given in Theorem  \ref{thmM2}. An auxiliary result in terms of a series associated with discrete matrix-valued shifts is proved, Lemma \ref{Blem1}. Theorem \ref{MainthmI} characterizes discrete matrix-valued frames of nonuniform shifts for the space $\ell^{2}\left(\Lambda,\mathcal{M}_n\right)$. Sec. 3 gives two perturbation results. Theorem \ref{PerMI} gives sufficient conditions for perturbation of  discrete matrix-valued frames of non-uniform shifts in the frequency domain. Relative error and frame conditions for a collection of  discrete matrix-valued  non-uniform shifts  are given in Theorem \ref{PerMIIx}. Some numerical examples are also given to illustrate our results.

\subsection{Basic Definitions and Notation}
In this section, we provide the necessary background on the structure of discrete matrix-valued nonuniform signal spaces along with some foundational results. As is conventional, the symbols $\mathbb{N}$, $\mathbb{Z}$, $\mathbb{R}$, and $\mathbb{C}$ represent the sets of natural numbers, integers, real numbers, and complex numbers, respectively. The set of all  $n\times n $ matrices having entries from complex field is  denoted by $\mathcal{M}_n$, it forms norm space with norm   $\|\mathbf{f}\|_{\mathcal{M}_{n}}=\sqrt{\sum_{m,k=1}^{n}\left|f_{m,k}\right|^{2}}$ for where $f_{m,k}$ denotes the $m^{th} $ row and $k^{th} $ column of vector $\mathbf{f}$. Throughout the paper  we use bold letters for matrices. Let $N \in \mathbb{N}$ and $r$ be an odd integer coprime with $N$, with the constraint $1 \leq r \leq 2 N-1$, assume $\Lambda=\left\{0, \frac{r}{N}\right\}+2 \mathbb{Z}$. For $p \in \mathbb{N}$, we write $\mathbb{N}_p = \{1,2,\dots, p\}$.

The following spaces will be used in the sequel: The space  $$\ell^{2}(\Lambda,\mathcal{M}_n):=\Bigg\{\mathbf{f}=
\big[f_{m, k}(\lambda)\big]_{ \lambda\in \Lambda \atop  1\le m,k \le n}: \sum\limits_{\substack{\lambda \in \Lambda \\ 1 \le m,k \le n}}|f_{m, k}(\lambda)|^2 < \infty \Bigg\},
$$
is a  Hilbert space with respect to the inner product
\[
\langle \mathbf{f}, \mathbf{g} \rangle = \sum_{\lambda \in \Lambda} \operatorname{tr} \left(\mathbf{f}(\lambda) \mathbf{g}(\lambda)^*  \right)= \sum\limits_{\substack{\lambda \in \Lambda \\ 1 \le m,k \le n}} f_{m,k}(\lambda)\overline{g_{m,k}(\lambda)}, \quad \mathbf{f}, \mathbf{g} \in \ell^{2}\left(\Lambda,\mathcal{M}_n\right),
\]
and its associated norm is
$$
\|\mathbf{f}\|_{\ell^{2}}=\sqrt{\sum_{\substack{\lambda \in \Lambda \\ 1 \le m,k \le n}}\left|{f}_{m,k}(\lambda)\right|^{2}}.
$$
It can be clearly seen that the above space will reduced to $\ell^2(\mathbb{Z}),$ if we take $N=1$ and $n=1$.

We recall a notation: For $p \in \mathbb{N}$, we denote $\mathbb{N}_p = \{1,2,\dots, p\}$. The space
	$$
	\ell^{2}\left(\Lambda, \mathbb{N}_p,\mathbb{C}\right):=\left\{\left\{c_{\lambda, j}\right\}_{\lambda\in \Lambda \atop j \in \mathbb{N}_p} \subset \mathbb{C} : \sum_{\substack{j =1 \\ \lambda \in \Lambda}}^{p}\left|c_{\lambda, j}\right|^{2}<\infty\right\}
	$$
	is a Hilbert space with respect to the following inner product
		\begin{align*}
		\left\langle\left\{c_{\lambda,  j}\right\}_{\lambda\in \Lambda \atop j \in \mathbb{N}_p},\left\{d_{\lambda, j}\right\}_{\lambda\in \Lambda \atop j \in \mathbb{N}_p}\right\rangle=\sum_{\substack{ j=1\\ \lambda \in \Lambda}}^{p} c_{\lambda, j} \overline{d_{\lambda,  j}}, \ \ \left\{c_{\lambda,  j}\right\}_{\lambda\in \Lambda \atop j \in \mathbb{N}_p},\left\{d_{\lambda, j}\right\}_{\lambda\in \Lambda \atop j \in \mathbb{N}_p}\in \ell^{2}\left(\Lambda,\mathbb{N}_p, \mathbb{C}\right),
	\end{align*}
	and the associated norm is given by
\begin{align*}
	\left\|\left\{c_{\lambda, j}\right\}_{\lambda\in \Lambda \atop j \in \mathbb{N}_p}\right\|_{\ell^{2}}=\sqrt{\sum_{\substack{ j=1\\ \lambda \in \Lambda}}^{p}\left|c_{\lambda, j}\right|^{2}}, \ \left\{c_{\lambda,  j}\right\}_{\lambda\in \Lambda \atop j \in \mathbb{N}_p} \in \ell^{2}\left(\Lambda, \mathbb{N}_p, \mathbb{C}\right).
\end{align*}
\begin{definition}
The \emph{displacement operator} or \emph{shift-operator} $\mathcal{D}_{2 N \lambda}: \ell^{2}\left(\Lambda,\mathcal{M}_n\right) \rightarrow \ell^{2}\left(\Lambda, \mathcal{M}_n\right)$,  for $\lambda \in \Lambda$, and  $\mathbf{f} \in \ell^{2}\left(\Lambda, \mathcal{M}_n\right)$,  is defined by
$$ \mathcal{D}_{2 N \lambda} \mathbf{f}=\left\{\left[\begin{array}{cccc}
	{f_{1,1}\left(\lambda^{\prime}-2 N \lambda\right)} &{f_{1,2}\left(\lambda^{\prime}-2 N \lambda\right)} & \cdots & {f_{1,n}\left(\lambda^{\prime}-2 N \lambda\right)}\\
	{f_{2,1}\left(\lambda^{\prime}-2 N \lambda\right)} &{f_{2,2}\left(\lambda^{\prime}-2 N \lambda\right)} & \cdots & {f_{2,n}\left(\lambda^{\prime}-2 N \lambda\right)}\\
	\vdots & \vdots & \vdots & \vdots \\
	{f_{n,1}\left(\lambda^{\prime}-2 N \lambda\right)} &{f_{n,2}\left(\lambda^{\prime}-2 N \lambda\right)} & \cdots & {f_{n,n}\left(\lambda^{\prime}-2 N \lambda\right)}
\end{array}\right]\right\}_{\lambda^{\prime} \in \Lambda}.
$$ This operator effectively shifts the argument of the function $\bold{f}$ by $2N\lambda$, mapping it within the same signal space.
\end{definition}
Let  $\left\{\mathbf{f}_{j}\right\}_{j=1}^{p} \subset \ell^{2}\left(\Lambda,\mathcal{M}_{n}\right)$  A collection of the  form $ \left\{ \mathcal{D}_{2 N \lambda} \mathbf{f}_{j}\right\}_{\lambda\in \Lambda \atop 1\le j \le p}$ is called a discrete \emph{matrix-valued nonuniform displacement  system} in the  discrete matrix-valued nonuniform sequence space $\ell^{2}\left(\Lambda, \mathcal{M}_n\right)$.

The  space
 $L^{2}\left(\Omega, \mathcal{M}_n\right):=\left\{\mathbf{F}=
\big[F_{m, k}(x)\big] : x\in \Omega, \,\, 1\le m,k \le n  :\int\limits_{\Omega}|F_{m, k}(x)|^2\,dx < \infty \right\}$  is a Hilbert space with respect to the inner product defined by 	
\begin{align*}
\langle\mathbf{F}, \mathbf{G}\rangle=\sum_{m,k=1}^{n} \int\limits_{\Omega}{{{F}_{m,k}(x)}}\overline{{G}_{m,k}(x)}\,  d x, \quad \mathbf{F}, \mathbf{G} \in L^{2}\left(\Omega, \mathcal{M}_n\right),
\end{align*}
 and the associated  norm is given by
\begin{align*}
\|\mathbf{F}\|_{L^{2}}=\sqrt{\sum_{m,k=1}^{n} \int_{\Omega}\left|{F}_{m,k}(x)\right|^{2}}, \ \mathbf{F} \in L^{2}\left(\Omega, \mathcal{M}_{n}\right).
\end{align*}

\begin{definition} The Fourier transform $\mathcal{F}: \ell^{2}\left(\Lambda, \mathcal{M}_n\right) \rightarrow L^{2}\left(\Omega, \mathcal{M}_n\right)$ is defined by
\begin{align*}
\mathcal{F}({\mathbf{f}})=\left[\begin{array}{cccc}
  	\sum\limits_{\lambda \in \Lambda}{f}_{1,1}(\lambda) e^{2 \pi i \lambda x}& \sum\limits_{\lambda \in \Lambda}{f}_{1,2}(\lambda) e^{2 \pi i \lambda x}&\cdots& \sum\limits_{\lambda \in \Lambda}{f}_{1,n}(\lambda) e^{2 \pi i \lambda x}\\
  	\sum\limits_{\lambda \in \Lambda}{f}_{2,1}(\lambda) e^{2 \pi i \lambda x}& \sum\limits_{\lambda \in \Lambda}{f}_{2,2}(\lambda) e^{2 \pi i \lambda x}&\cdots& \sum\limits_{\lambda \in \Lambda}{f}_{2,n}(\lambda) e^{2 \pi i \lambda x}\\
  	\vdots&\vdots&\vdots&\vdots\\
  	\sum\limits_{\lambda \in \Lambda}{f}_{n,1}(\lambda) e^{2 \pi i \lambda x}& \sum\limits_{\lambda \in \Lambda}{f}_{n,2}(\lambda) e^{2 \pi i \lambda x}&\cdots& \sum\limits_{\lambda \in \Lambda}{f}_{n,n}(\lambda) e^{2 \pi i \lambda x}
  \end{array}\right].
\end{align*}
It is well defined, since $(\Lambda, \Omega)$ is a spectral pair by Example \ref{exp12}, so series corresponding to each coordinate will converges in  $L^{2}(\Omega, \mathcal{M}_n)$.
  \end{definition}
  Moreover,
  $\mathcal{F}$
   is bijective, with the inverse given by the inverse Fourier transform $\mathcal{F}^{-1} : L^{2}\left(\Omega, \mathcal{M}_{n}\right) \rightarrow \ell^{2}\left(\Lambda,\mathcal{M}_{n}\right)$ is given by

  $$
  \mathcal{F}^{-1}({\mathbf{F}})=\left[\begin{array}{cccc}
  	\int_{\Omega}{F}_{1,1}(x) e^{-2 \pi i \lambda x}\, dx & \int_{\Omega}{F}_{1,2}(x) e^{-2 \pi i \lambda x}dx &\cdots& \int_{\Omega}{F}_{1,n}(x) e^{-2 \pi i \lambda x}dx\\
  	\int_{\Omega}{F}_{2,1}(x) e^{-2 \pi i \lambda x}\, dx & \int_{\Omega}{F}_{2,2}(x) e^{-2 \pi i \lambda x}dx &\cdots& \int_{\Omega}{F}_{2,n}(x) e^{-2 \pi i \lambda x}dx\\
  	\vdots&\vdots&\vdots&\vdots\\
  	\int_{\Omega}{F}_{n,1}(x) e^{-2 \pi i \lambda x}\, dx & \int_{\Omega}{F}_{n,2}(x) e^{-2 \pi i \lambda x}dx &\cdots& \int_{\Omega}{F}_{n,n}(x) e^{-2 \pi i \lambda x}dx
  \end{array}\right]
  $$

  Now  we provide some  identities that are relevant  to the Fourier transform. Let $\mathbf{f}, \mathbf{g} \in \ell^{2}\left(\Lambda,\mathcal{M}_n\right)$ and $\lambda \in \Lambda$. Then,
\begin{enumerate}
\item $\mathcal{F}\left(\mathcal{D}_{2 N \lambda} \mathbf{f}\right)(x)=e^{4 \pi i N \lambda x} \mathcal{F}(\mathbf{f})(x)$.

\item  $\langle\mathbf{f}, \mathbf{g}\rangle_{\ell^{2}}=\langle\mathcal{F}(\mathbf{f}), \mathcal{F}(\mathbf{g})\rangle_{L^{2}}$ (Parseval's relation).

\item $\|\mathbf{f}\|_{\ell^{2}}=\|\mathcal{F}(\mathbf{f})\|_{L^{2}}$ (Plancherel's formula).
\end{enumerate}

\section{Matrix-valued discrete frames of  non-uniform shifts}
In this section, we study frame properties for a collection of non-uniform shifts of the form $ \left\{ \mathcal{D}_{2 N \lambda} \mathbf{f}_{j}\right\}_{\lambda\in \Lambda \atop j \in \mathbb{N}_p}$ in the space $\ell^{2}\left(\Lambda,\mathcal{M}_n\right)$. We begin with the following definition:
\begin{definition}\label{defnM1}
For each $j \in \mathbb{N}_p$,  let $\mathbf{f}_{j}$ be a  non-zero element of the space $\ell^{2}\left(\Lambda,\mathcal{M}_{n}\right)$.
A collection  $ \left\{ \mathcal{D}_{2 N \lambda} \mathbf{f}_{j}\right\}_{\lambda\in \Lambda \atop j \in \mathbb{N}_p} \subset \ell^{2}\left(\Lambda,\mathcal{M}_n\right)$ is  a \emph{matrix-valued discrete frame of  non-uniform shifts}  (\emph{matrix-valued DFNS}, in short) for $\ell^{2}\left(\Lambda,\mathcal{M}_n\right)$,  if there exist  $a_0$, $b_0 \in (0, \infty)$ such that
	\begin{align}\label{defn1.1}
		a_0\|\mathbf{f}\|_{\ell^{2}}^{2} \leq \sum_{\substack{ j=1 \\ \lambda \in \Lambda}}^{p}\left|\left\langle\mathbf{f},\mathcal{D}_{2 N \lambda} \mathbf{f}_{j}\right\rangle\right|^{2} \leq b_0\|\mathbf{f}\|_{\ell^{2}}^{2} \text { for all } \mathbf{f} \in \ell^{2}\left(\Lambda,\mathcal{M} _{n}\right).
	\end{align}
\end{definition}
Positive scalars $a_0$ and $b_0$ in the \emph{frame inequality} \eqref{defn1.1}  are known as the \emph{lower frame bound} and \emph{upper frame bound} of the frame $ \left\{ \mathcal{D}_{2 N \lambda} \mathbf{f}_{j}\right\}_{\lambda\in \Lambda \atop j \in \mathbb{N}_p}$, respectively. If only the right-hand side inequality in \eqref{defn1.1} is satisfied, then we say that $ \left\{ \mathcal{D}_{2 N \lambda} \mathbf{f}_{j}\right\}_{\lambda\in \Lambda \atop j \in \mathbb{N}_p}$ is a \emph{ matrix-valued discrete Bessel sequence of non-uniform shifts} with  \emph{Bessel bound} $b_0$. If $ \left\{ \mathcal{D}_{2 N \lambda} \mathbf{f}_{j}\right\}_{\lambda\in \Lambda \atop j \in \mathbb{N}_p}$ matrix-valued discrete Bessel sequence of non-uniform shifts, then
the map $\mathcal{T}: \ell^{2}\left(\Lambda, \mathbb{N}_p,\mathbb{C}\right) \rightarrow \ell^{2}\left(\Lambda, \mathcal{M}_{n}\right)$ defined by
	$$\mathcal{T}\left(\left\{c_{\lambda, j}\right\}_{\lambda\in \Lambda \atop j \in \mathbb{N}_p}\right)=\sum\limits_{\substack{ j=1\\ \lambda \in \Lambda}}^{p} c_{\lambda, j}\mathcal{D}_{2 N \lambda} \mathbf{f}_{j},$$
	is called the  \emph{pre-frame operator} (or \emph{synthesis operator}) of  $ \left\{ \mathcal{D}_{2 N \lambda} \mathbf{f}_{j}\right\}_{\lambda\in \Lambda \atop j \in \mathbb{N}_p}$ and the Hilbert-adjoint operator of $\mathcal{T}$ is the map $\mathcal{T}^{*}$ : $\ell^{2}\left(\Lambda,\mathcal{M}_{n}\right) \rightarrow \ell^{2}\left(\Lambda, \mathbb{N}_p, \mathbb{C}\right)$  given by
	$$\mathcal{T}^{*}(\mathbf{f})= \left\{\left\langle\mathbf{f}, \mathcal{D}_{2 N \lambda} \mathbf{f}_{j}\right\rangle\right\}_{\lambda\in \Lambda \atop j \in \mathbb{N}_p}. $$
$\mathcal{T}^{*}$ is  called the \emph{analysis operator} of $ \left\{ \mathcal{D}_{2 N \lambda} \mathbf{f}_{j}\right\}_{\lambda\in \Lambda \atop j \in \mathbb{N}_p}$ .
	The composition $\mathcal{S}=\mathcal{T}\mathcal{T}^{*}: \ell^{2}\left(\Lambda,\mathcal{M}_{n}\right) \rightarrow \ell^{2}\left(\Lambda,\mathcal{M}_{n}\right)$ is called the  \emph{frame operator}, which is given by
\begin{align*}
	\mathcal{S} =\sum_{\substack{j=1,\\ \lambda \in \Lambda}}^p\left\langle\mathbf{f},\mathcal{D}_{2 N \lambda} \mathbf{f}_{j}\right\rangle \mathcal{D}_{2 N \lambda} \mathbf{f}_{j}.
\end{align*}
\begin{proposition} [\cite{Heil2}, Theorem 7.4] \label{Bes1}Let $\left\{g_{k}\right\}_{k \in \mathbb{I}}$ be a sequence of vectors in $\mathcal{H}$, such that there exists a constant $b_{0}>0$ satisfying
\begin{align*}
	\left\|\sum\limits_{k \in \mathbb{\mathbb{I}}} a_{k} g_{k}\right\|_{\mathcal{H}}^{2} \leq b_{0} \sum\left|a_{k}\right|^{2}
\end{align*}
	for all finite sequences $\left\{a_{k}\right\}$, then the series $\sum_{k \in \mathbb{I}} a_{k} g_{k}$ converges for every $\left\{a_{k}\right\}_{k \in \mathbb{I}} \in$ $\ell^{2}(\mathbb{I})$, further $\left\{g_{k}\right\}_{k \in \mathbb{I}}$ is a Bessel sequence with Bessel bound $b_{0}$.\\
\end{proposition}
The following lemma  will be used in the sequel. It easily  follows from the mathematical induction method.
 \begin{lemma}\label{Blem1} $\left|\sum\limits_{k=1}^{t} w_{k}\right|^{2} \leq 2^{t-1} \sum\limits_{k=1}^{t}\left|w_{k}\right|^{2}$, where $w_k \in \mathbb{C}$, $1\leq  k \leq t$.
	\end{lemma}
Our first result provides a sufficient condition for matrix-valued discrete Bessel sequence of  non-uniform shifts in terms of the Fourier transforms of signals.
\begin{theorem}\label{thmM1}
For each $j \in \mathbb{N}_p$, let  $\mathbf{f}_{j}$ be a non-zero element in the space $\ell^{2}\left(\Lambda,\mathcal{M}_{n}\right)$, and let  $b_o \in (0, \infty)$ be  such that
\begin{align*}
	\left\|\mathcal{F}\left(\mathbf{f}_{j}(x)\right)\right\|_{\mathcal{M}_{n}} \leq b_o, \  \text{a.e.} \   x \in \Omega, \,   j \in \mathbb{N}_p.
\end{align*}
	Then, the collection $\left\{ \mathcal{D}_{2 N \lambda} \mathbf{f}_{j} \right\}_{\lambda\in \Lambda \atop j \in \mathbb{N}_p}$ form a matrix-valued discrete Bessel sequence of  non-uniform shifts with Bessel bound $2^{p-1}b_{o}^2n^2.$
\end{theorem}
\proof For any finite sequence  $\left\{c_{\lambda,  j}\right\}_{\lambda\in \Lambda \atop j \in \mathbb{N}_p} \in \ell^{2}\left(\Lambda,
\left[p\right] ,\mathbb{C}\right)$, we have
	\begin{align}\label{eq3.4}
		\left\|\sum_{\substack{j =1 \\ \lambda \in \Lambda}}^{p} c_{\lambda,j}\mathcal{D}_{2N\lambda}\mathbf{f}_{j}\right\|_{\ell^{2}}
		&= \left\|\mathcal{F}\left(\sum_{\substack{j =1 \\ \lambda \in \Lambda}}^{p} c_{\lambda,j} \mathcal{D}_{2 N \lambda} \mathbf{f}_{j}\right)\right\|_{L^{2}} \notag \\
		&\leq \sum_{j=1}^p \left\|\sum_{\lambda \in \Lambda} c_{\lambda, j} e^{4 \pi i N \lambda x} \mathcal{F}\left(\mathbf{f}_{j}\right)(x)\right\|_{L^{2}} \notag \\
		&= n b_{o}\sum_{j=1}^p \sqrt{\int_{\Omega} \left| \sum_{\lambda \in \Lambda} c_{\lambda, j} e^{4 \pi i N \lambda x} \right|^{2} d x}.
	\end{align}
	Using \eqref{eq3.4} and the fact  that   $\left\{e^{4 \pi i N \lambda x}\right\}_{\lambda \in \Lambda}$ is orthonormal system in $L^{2}(\Omega)$, we get
	\begin{align}\label{eqd1}
		\left\|\sum_{\substack{j =1 \\ \lambda \in \Lambda}}^{p} c_{\lambda, j} \mathcal{D}_{2 N \lambda} \mathbf{f}_{j}\right\|_{\ell^{2}} \leq n b_{o} \sum_{j=1}^p \sqrt{\sum_{\lambda \in \Lambda} \left|c_{\lambda,j}\right|^{2}}
	\end{align}
	Therefore, by inequality  \eqref{eqd1} and Lemma \ref{Blem1}, we get	
	\begin{align}\label{xxrr1}
		\left\|\sum_{\substack{j =1 \\ \lambda \in \Lambda}}^{p}c_{\lambda, j} \mathcal{D}_{2 N \lambda} \mathbf{f}_{j}\right\|_{\ell^{2}}^{2}
		\leq \left( n b_{o}\sum_{\substack{j=1}}^p \sqrt{\sum_{\lambda \in \Lambda} \left|c_{\lambda, j}\right|^{2}}\right)^{2}
		\leq  b_{o}^{2}2^{p-1}  n^{2} \sum_{\substack{j=1}}^p \sum_{\lambda \in \Lambda} \left|c_{\lambda, j}\right|^{2}.
	\end{align}	
Inequality \eqref{xxrr1} and Proposition \ref{Bes1} gives the result with the required Bessel bound.
	\endproof
The following example demonstrates Theorem \ref{thmM1}.
\begin{example}\label{Exam1}
 Let $\Lambda=\left\{0, \frac{1}{2}\right\}+2 \mathbb{Z}, \Omega=$ $\left[0, \frac{1}{2}\right)\cup\left[1, \frac{3}{2}\right).$ Define $\left\{\mathbf{f}_{j}\right\}_{j=1}^{8} \subset \ell^{2}\left(\Lambda,\mathcal{M}_{2}\right)$   as below:
\begin{align*}
		& \mathbf{f}_{1}(0)=\left[\begin{array}{cc}
			1 &0\\
			0&1
		\end{array}\right], \quad \mathbf{f}_{1}(4)=\left[\begin{array}{cc}
			0&1\\
			1&0
		\end{array}\right], \quad \mathbf{f}_{1}(\lambda)=\left[\begin{array}{cc}
			0 & 0\\
			0&0
		\end{array}\right] \text { for } \lambda \in \Lambda \backslash\{0,4\} ;\\
		& \mathbf{f}_{2}(0)=\left[\begin{array}{cc}
			1&0\\
			0&-1
		\end{array}\right], \quad \mathbf{f}_{2}(4)=\left[\begin{array}{cc}
			0&-i\\
			i&0
		\end{array}\right], \quad \mathbf{f}_{2}(\lambda)=\left[\begin{array}{cc}
			0 &0\\
			0&0
		\end{array}\right] \text { for } \lambda \in \Lambda \backslash\{0,4\};\\
		&\mathbf{f}_{3}(0)=\left[\begin{array}{cc}
			0 & i\\
			-i&0
		\end{array}\right], \quad \mathbf{f}_{3}(4)=\left[\begin{array}{cc}
			0 & 1\\
			1 & 0
		\end{array}\right], \quad \mathbf{f}_{3}(\lambda)=\left[\begin{array}{cc}
			0 & 0\\
			0& 0
		\end{array}\right] \text { for } \lambda \in \Lambda \backslash\{0,4\} ;\\
		&\mathbf{f}_{4}(0)=\left[\begin{array}{cc}
			0 & 1\\
			1&0
		\end{array}\right], \quad \mathbf{f}_{4}(4)=\left[\begin{array}{cc}
			1 & 0\\
			0&-1
		\end{array}\right], \quad \mathbf{f}_{4}(\lambda)=\left[\begin{array}{cc}
			0 & 0\\
			0& 0
		\end{array}\right] \text { for } \lambda \in \Lambda \backslash\{0,4\}; \\
		& \mathbf{f}_{5}(\frac{1}{2})=\left[\begin{array}{cc}
			1 & 0\\
			0&1
		\end{array}\right], \quad \mathbf{f}_{5}(\frac{1}{2}+4)=\left[\begin{array}{cc}
			0 & 1\\
			1 &0
		\end{array}\right], \quad \mathbf{f}_{5}(\lambda)=\left[\begin{array}{cc}
			0 & 0\\
			0& 0
		\end{array}\right] \text { for } \lambda \in \Lambda \backslash\{\frac{1}{2},\frac{1}{2}+4\} ;\\
		& \mathbf{f}_{6}(\frac{1}{2})=\left[\begin{array}{cc}
			1 & 0\\
			0&-1
		\end{array}\right], \quad \mathbf{f}_{6}(\frac{1}{2}+4)=\left[\begin{array}{cc}
			0 & -i\\
			i &0
		\end{array}\right], \quad \mathbf{f}_{6}(\lambda)=\left[\begin{array}{cc}
			0 & 0\\
			0& 0
		\end{array}\right] \text { for } \lambda \in \Lambda \backslash\{\frac{1}{2},\frac{1}{2}+4\}; \\
		& \mathbf{f}_{7}(\frac{1}{2})=\left[\begin{array}{cc}
			0 & i\\
			-i&0
		\end{array}\right], \quad \mathbf{f}_{7}(\frac{1}{2}+4)=\left[\begin{array}{cc}
			0 & 1\\
			1 &0
		\end{array}\right], \quad \mathbf{f}_{7}(\lambda)=\left[\begin{array}{cc}
			0 & 0\\
			0& 0
		\end{array}\right] \text { for } \lambda \in \Lambda \backslash\{\frac{1}{2},\frac{1}{2}+4\} ;\\
		& \mathbf{f}_{8}(\frac{1}{2})=\left[\begin{array}{cc}
			0& 1\\
			1&0
		\end{array}\right], \quad \mathbf{f}_{8}(\frac{1}{2}+4)=\left[\begin{array}{cc}
			1 & 0\\
			0 &-1
		\end{array}\right], \quad \mathbf{f}_{8}(\lambda)=\left[\begin{array}{cc}
			0 & 0\\
			0& 0
		\end{array}\right] \text { for } \lambda \in \Lambda \backslash\{\frac{1}{2},\frac{1}{2}+4\}; \\
\end{align*}
	Now for $1\le j \le 8$, the corresponding Fourier transforms $\mathcal{F}\left(\mathbf{f}_{j}\right)(x)$ are given by
	\begin{equation*}
		\begin{aligned}
			\mathcal{F}\left( \mathbf{f}_{1}\right)(x)&=\left[\begin{array}{@{}cc@{}}
				1 & e^{8\pi i x} \\
				e^{8\pi i x} & 1
			\end{array}\right], &
			\mathcal{F}\left(\mathbf{f}_{2}\right)(x)&=\left[\begin{array}{@{}cc@{}}
				1 & -ie^{8\pi i x} \\
				ie^{8\pi i x} & -1
			\end{array}\right], \\
			\mathcal{F}\left(\mathbf{f}_{3}\right)(x)&=\left[\begin{array}{@{}cc@{}}
				0 & i+e^{8 \pi i x}\\
				-i+e^{8 \pi i x} & 0
			\end{array}\right], &
			\mathcal{F}\left(\mathbf{f}_{4}\right)(x)&=\left[\begin{array}{@{}cc@{}}
				e^{8 \pi i x} & 1 \\
				1  & -e^{8 \pi i x}
			\end{array}\right], \\
			\mathcal{F}\left(\mathbf{f}_{5}\right)(x)&=\left[\begin{array}{@{}cc@{}}
				e^{\pi i x} &  e^{9\pi i x} \\
				e^{9\pi i x} & e^{\pi i x}
			\end{array}\right], &
			\mathcal{F}\left(\mathbf{f}_{6}\right)(x)&=\left[\begin{array}{@{}cc@{}}
				e^{\pi i x} &  - ie^{9\pi i x} \\
				ie^{9\pi i x} & -e^{\pi i x}
			\end{array}\right], \\
			\mathcal{F}\left(\mathbf{f}_{7}\right)(x)&=\left[\begin{array}{@{}cc@{}}
				0 & ie^{\pi i x}+e^{9\pi i x}\\
				-ie^{\pi i x} + e^{9\pi i x} & 0
			\end{array}\right], &
			\mathcal{F}\left( \mathbf{f}_{8}\right)(x)&=\left[\begin{array}{@{}cc@{}}
				e^{9\pi i x}& e^{\pi i x}\\
				e^{\pi i x}  & - e^{9\pi i x}
			\end{array}\right].
		\end{aligned}
	\end{equation*}
	By using Lemma \ref{Blem1}, it is easy to calculate for $ 1\le j \le 8$ that
	$\left\|\mathcal{F}\left(\mathbf{f}_{j}\right)(x)\right\|_{\mathcal{M}_{2} }\leq 2$.
	Thus, by Theorem \ref{thmM1}, the collection of matrix-valued sequences  $\left\{ \mathcal{D}_{2 N \lambda} \mathbf{f}_{j}\right\}_{\lambda\in \Lambda \atop j \in \mathbb{N}_8}$
 is a Bessel sequence with Bessel bound $2^{11}$.
\end{example}
\begin{theorem}[\cite{Heil2}, Theorem 7.4]\label{RefH1} A sequence $\left\{g_{k}\right\}_{k \in \mathbb{I}}$ of vectors in a Hilbert space $\mathcal{H}$ is a Bessel sequence with Bessel bound $b_{0}$ if and only if the pre-frame operator $\mathcal{T}: \ell^{2}(\mathbb{I}) \rightarrow \mathcal{H}$, $\mathcal{T}\left(\left\{a_{k}\right\}_{k \in \mathbb{I}}\right)=\sum\limits_{k \in \mathbb{I}} a_{k} g_{k}$ is bounded operator from $\ell^{2}(\mathbb{I})$ into $\mathcal{H}$ with $\|\mathcal{T}\|^{2} \leq b_{0}$.
\end{theorem}

\begin{lemma}\cite{Heil2}\label{RefH2} For any $b>0$ and any constant $c_{0} \in \mathbb{R},\left\{\sqrt{b} e^{2 \pi i b l\left(\xi+c_{0}\right)}\right\}_{l \in \mathbb{Z}}$ is a complete orthonormal system of $L^{2}\left(0, \frac{1}{b}\right)$.
\end{lemma}
Theorem \ref{RefH1}  and Lemma \ref{RefH2} are used in the following necessary condition for the existence of matrix-valued discrete Bessel sequence of  non-uniform shifts.
\begin{theorem}\label{thmM2}
Let $\left\{ \mathcal{D}_{2 N \lambda} \mathbf{f}_{j}\right\}_{\lambda\in \Lambda \atop j \in \mathbb{N}_p}$ be a matrix-valued discrete Bessel sequence of  non-uniform shifts in  $\ell^2(\Lambda,\mathcal{M}_n)$  with Bessel bound $b_{0}$.
Then, for a.e. $x \in \Omega$ and $j \in \mathbb{N}_p$, we have
\begin{align*}
\left\|\mathcal{F}\left( \mathbf{f}_{j}\right)(x)\right\|_{\mathcal{M}_{n}} \leq N + b_0.
\end{align*}
\end{theorem}
\proof  Using  Lemma \ref{RefH2}, for any $h \in L^{2}\left(0, \frac{1}{4 N}\right)$,
there exists sequence $\{d_{l}\}_{l \in \mathbb{Z}}$ of scalars  such that
	$h(x)=\sum_{l \in \mathbb{Z}} d_{l} \sqrt{4 N} e^{2 \pi i(4 N) lx} $
	with $\int_{0}^{\frac{1}{4 N}}|h(x)|^{2} d x=\sum_{l \in \mathbb{Z}}\left|d_{l}\right|^{2}<\infty$.
	Let $\bold{f}_t\in \{\bold{f}_1,\bold{f}_2,\dots,\bold{f}_p\}$ be an arbitrary element.
	Choose a complex  sequence $\left\{c_{\lambda, j}\right\}_{\lambda\in \Lambda \atop j \in \mathbb{N}_p}$ 	as  below
\begin{align}\label{caser1}
	c_{\lambda,  j}= \begin{cases}d_{l}, & \text { if } \lambda=2 l, l \in \mathbb{Z},  j=t, \\ 0, & \text { otherwise }.
\end{cases}
\end{align}
	Clearly  $\left\{c_{\lambda, j}\right\}_{\lambda\in \Lambda \atop j \in \mathbb{N}_p} \in \ell^{2}\left(\Lambda, \mathbb{N}_p, \mathbb{C}\right)$.  As $	\left\{ \mathcal{D}_{2 N \lambda} \mathbf{f}_{j}\right\}_{\lambda\in \Lambda \atop j \in \mathbb{N}_p} $ is Bessel sequence with Bessel bound $b_{0}$, by Theorem \ref{RefH1}, the associated pre-frame operator $\mathcal{T}: \ell^{2}\left(\Lambda, \mathbb{N}_p, \mathbb{C}\right)
	\rightarrow \ell^{2}\left(\Lambda,\mathcal{M}_{n}\right)$ is bounded and satisfies
	\begin{align}\label{caser2}
		\left\|\mathcal{T}\left\{c_{\lambda,j}\right\}\right\|_{\ell^{2}}^{2} \leq b_{0} \sum_{\substack{ j=1\\ \lambda \in \Lambda}}^p\left|c_{\lambda,  j}\right|^{2}.
\end{align}
	Using \eqref{caser1}, we compute
	\begin{align}\label{caser3}
		\left\|\mathcal{T}\left\{c_{\lambda, j}\right\}\right\|_{\ell^{2}}^{2}= & \left\|\sum_{l \in \mathbb{Z}} d_{l} \mathcal{D}_{4 N l} \mathbf{f}_{t}\right\|_{\ell^{2}}^{2} \notag\\
		= & \left\|\mathcal{F}\left(\sum_{l \in \mathbb{Z}} d_{l} \mathcal{D}_{4 N l} \mathbf{f}_{{t}}\right)\right\|_{L^{2}}^{2}  \notag\\
		= & \left\|\sum_{l \in \mathbb{Z}} d_{l} e^{4 \pi i N(2 l)\left(x\right)} \mathcal{F}\left( \mathbf{f}_{{t}}\right)(x)\right\|_{L^{2}}^{2} \notag\\
		= & \int_{\Omega} \sum_{m,k=1}^{n}\left|\sum_{l \in \mathbb{Z}} d_{l} e^{4 \pi i N(2 l)\left(x\right)}\mathcal{F}\left( \mathbf{f}_{t}\right)_{m,k}(x)\right|^{2} dx \notag\\
		= & \sum_{m,k=1}^{n} \int_{0}^{\frac{1}{2}}\left|\sum_{l \in \mathbb{Z}} d_{l} e^{4 \pi i N(2 l)\left(x\right)}\right|^{2}\left(\left|\mathcal{F}\left( \mathbf{f}_{t}\right)_{m,k}(x)\right|^{2} +\left|\mathcal{F}\left( \mathbf{f}_{t}\right)_{m,k}\left(x+\frac{N}{2}\right)\right|^{2}\right) dx \notag\\
		= & \sum_{m,k=1}^{n} \sum_{g=0}^{2 N-1} \int_{0}^{\frac{1}{4 N}}\left|\sum_{l \in \mathbb{Z}} d_{l} e^{4 \pi i N(2 l)\left(x\right)}\right|^{2}\left(\left|\mathcal{F}\left( \mathbf{f}_{t}\right)_{m,k}\left(x+\frac{g}{4 N}\right)\right|^{2}\right. \notag\\
		& \left.+\left|\mathcal{F}\left(\mathbf{f}_{t}\right)_{m,k}\left(x+\frac{N}{2}+\frac{g}{4 N}\right)\right|^{2}\right) d x.
	\end{align}
	Utilizing equations \eqref{caser1}, \eqref{caser2}and \eqref{caser3}, we obtain.
	\begin{align*}
 &\sum_{m,k=1}^{n} \sum_{g=0}^{2 N-1} \int_{0}^{\frac{1}{4 N}}|h(x)|^{2}\left(\left|\mathcal{F}\left(\mathbf{f}_{t}\right)_{m,k}\left(x+\frac{g}{4 N}\right)\right|^{2}+\left|\mathcal{F}\left(\mathbf{f}_{t}\right)_{m,k}\left(x+\frac{N}{2}+\frac{g}{4 N}\right)\right|^{2}\right) dx\\
&\leq 4 N b_{0}\int_{0}^{\frac{1}{4 N}}|h(\xi)|^{2} d \xi.
\end{align*}
	As $h$ is an arbitrary element in the space $L^{2}\left(0, \frac{1}{4 N}\right)$, it follows that for a.e. $ x\in\left[0, \frac{1}{4 N}\right)$,
	$$
	\sum_{m,k=1}^{n} \sum_{g=0}^{2 N-1}\left|\mathcal{F}\left( \mathbf{f}_{t}\right)_{m,k}\left(x+\frac{g}{4 N}\right)\right|^{2}+\left|\mathcal{F}\left(\mathbf{f}_{t}\right)_{m,k}\left(x+\frac{N}{2}+\frac{g}{4 N}\right)\right|^{2} \leq 4 N b_{0},
	$$
	which implies
\begin{align*}
	\sum_{m,k=1}^{n}\left|\mathcal{F}\left(\mathbf{f}_{j}\right)_{m,k}(x)\right|^{2} \leq 4 N b_{0} \  \text { for a.e. } x \in \Omega, \,  j \in \mathbb{N}_p.
\end{align*}
This gives $\left\|\mathcal{F}\left( \mathbf{f}_{j}\right)(x)\right\|_{\mathcal{M}_{n}} \leq 2 \sqrt{N b_0} \leq N + b_0$ for a.e. $x \in \Omega$, $j \in \mathbb{N}_p$.
	This completes the proof.
\endproof
\begin{corollary}\label{corM3}
 $\left\{ \mathcal{D}_{2 N \lambda} \mathbf{f}_{j} \right\}_{\lambda\in \Lambda \atop j \in \mathbb{N}_p} \subset  \ell^2(\Lambda,\mathcal{M}_n)$ is a matrix-valued discrete Bessel sequence of non-uniform shifts if and only if there exists there exists a positive finite constant $b_{0}$ such that
\begin{align}\label{eqcor3}
		\left\|\mathcal{F}\left(\mathbf{f}_{j}(x)\right)\right\|_{\mathcal{M}_{n}} \leq  b_{0}, \  \text{a.e. } x \in \Omega, \,  j \in \mathbb{N}_p.
\end{align}
\end{corollary}
\begin{remark} Condition \eqref{eqcor3} does not guarantee the lower frame condition for system   $\left\{ \mathcal{D}_{2 N \lambda} \mathbf{f}_{j}\right\}_{\lambda\in \Lambda \atop j \in [p]}$. Indeed,  consider $\Lambda=\left\{0, \frac{r}{N}\right\}+2 \mathbb{Z}$, $\Omega=\left[0, \frac{1}{2}\right)\cup\left[\frac{N}{2}, \frac{N+1}{2}\right),p=2 $ and  define $\left\{\mathbf{f}_{j}\right\}_{j=1}^{2} \subset \ell^{2}\left(\Lambda,\mathcal{M}_{2}\right)$  in term of corresponding Fourier transforms $\mathcal{F}\left(\mathbf{f}_{j}\right)(x)$ as follows:
\begin{align*}
	\mathcal{F}\left(\mathbf{f}_{1}\right)=\left[\begin{array}{cc}
		\sqrt{2 N} \phi_{\left[0, \frac{1}{4N}\right)}(x) & 0 \\
		0 &\sqrt{2 N} \phi_{\left[0, \frac{1}{4N}\right)}(x)
	\end{array}\right] \text {and } \mathcal{F}\left(\mathbf{f}_{2}\right)=\left[\begin{array}{cc}
		0 & \sqrt{2 N} \phi_{\left[0, \frac{1}{4 N}\right)}(x)\\
		\sqrt{2 N} \phi_{\left[0, \frac{1}{4 N}\right)}(x) &0
	\end{array}\right].
\end{align*}
A simple computation gives $\left\|\mathcal{F}\left( \mathbf{f}_{j}\right)(x)\right\|_{\mathcal{M}_{2}} \leq 2\sqrt{ N}<\infty$  for a.e. $x \in \Omega$ and $1\le j \le 2 $. Thus, by Corollary \ref{corM3}, the  system of matrix-valued discrete non-uniform shifts
$\left\{ \mathcal{D}_{2 N \lambda} \mathbf{f}_{j}\right\}_{\lambda\in \Lambda \atop j \in [p]}$  is a Bessel sequence in $\ell^{2}\left(\Lambda, \mathcal{M}_{2}\right)$.
Suppose $\left\{ \mathcal{D}_{2 N \lambda} \mathbf{f}_{j}\right\}_{\lambda\in \Lambda \atop j \in [p]}$ satisfies the lower frame condition in the space $\ell^{2}\left(\Lambda, \mathcal{M}_{2}\right)$.
Then, there exists a positive real number $a_{0}$  such that
\begin{align}\label{euu12}
a_{0}\|\mathbf{f}\|^{2} \leq \sum_{\substack{ j=1 \\ \lambda \in \Lambda}}^{2}\left|\left\langle\mathbf{f},  \mathcal{D}_{2 N \lambda} \mathbf{f}_{j}\right\rangle\right|^{2} \ \text{for all} \ \mathbf{f} \in \ell^{2}\left(\Lambda,\mathcal{M}_{2}\right).
	\end{align}
For any  $\mathbf{f} \in \ell^{2}\left(\Lambda,\mathcal{M}_{2}\right)$, we compute
	\begin{align}\label{eqotto1}
		& \sum_{\substack{ j=1\\
				\lambda \in \Lambda}}^2\left|\left\langle\mathbf{f}, \mathcal{D}_{2 N \lambda} \mathbf{f}_{j}\right\rangle\right|^{2} \notag\\
		= & \sum_{\substack{j=1 \\ \lambda \in\Lambda}}^2 \left|\left\langle\mathcal{F}(\mathbf{f}),\mathcal{F}\left(\mathcal{D}_{2 N \lambda} \mathbf{f}_{j}\right)\right\rangle\right|^{2} \notag \\
		= & \sum_{\substack{j=1 \\  \lambda \in \Lambda}}^2\left|\left\langle\mathcal{F}(\mathbf{f})(x), e^{4 \pi i N \lambda x} \mathcal{F}\left(\mathbf{f}_{j}\right)(x)\right\rangle\right|^{2} \notag\\
		= & \sum_{\lambda \in \Lambda} \left\lvert\, \int_{\Omega}\mathcal{F}(\mathbf{f})_{1,1}(x) e^{-4 \pi i N \lambda x} \sqrt{2 N} \phi_{\left[0, \frac{1}{4 N}\right[}(x) dx\right|^{2} \notag\\
		& +\sum_{\lambda \in \Lambda}\left|\int_{\Omega}\mathcal{F}(\mathbf{f})_{1,2}(x) e^{-4 \pi i N \lambda x} \sqrt{2 N} \phi_{\left[0, \frac{1}{4 N}\right[}(x) dx\right|^{2} \notag\\
		& +\sum_{\lambda \in \Lambda}\left|\int_{\Omega}\mathcal{F}(\mathbf{f})_{2,1}(x) e^{-4 \pi i N \lambda x} \sqrt{2 N} \phi_{\left[0, \frac{1}{4 N}\right[}(x) dx\right|^{2} \notag\\
		& +\sum_{\lambda \in \Lambda}\left|\int_{\Omega}\mathcal{F}(\mathbf{f})_{2,2}(x) e^{-4 \pi i N \lambda x} \sqrt{2 N} \phi_{\left[0, \frac{1}{4 N}\right[}(x) dx\right|^{2} \notag\\
		= & 2 N\left(\sum_{l\in\mathbb{Z}}\left|\int_{0}^{\frac{1}{4 N}}\mathcal{F}(\mathbf{f})_{1,1}(x) e^{-2 \pi i(4 N) l x} dx\right|^{2}+\sum_{l \in \mathbb{Z}}\left|\int_{0}^{\frac{1}{4 N}}\mathcal{F}(\mathbf{f})_{1,1}(x) e^{-4 \pi i r x} e^{-2 \pi i(4 N) l x} dx\right|^{2}\right. \notag\\
		& \left.+\sum_{l \in \mathbb{Z}}\left|\int_{0}^{\frac{1}{4 N}}\mathcal{F}(\mathbf{f})_{1,2}(x) e^{-2 \pi i(4 N) lx} dx\right|^{2}+\sum_{l \in \mathbb{Z}}\left|\int_{0}^{\frac{1}{4 N}}\mathcal{F}(\mathbf{f})_{1,2}(x) e^{-4 \pi i rx} e^{-2 \pi i(4 N) lx} dx\right|^{2}\right)\notag\\
		& \left.+\sum_{l \in \mathbb{Z}}\left|\int_{0}^{\frac{1}{4 N}}\mathcal{F}(\mathbf{f})_{2,1}(x) e^{-2 \pi i(4 N) l x} d x\right|^{2}+\sum_{z \in \mathbb{Z}}\left|\int_{0}^{\frac{1}{4 N}}\mathcal{F}(\mathbf{f})_{2,1}(x) e^{-4 \pi i r x} e^{-2 \pi i(4 N) lx} d x\right|^{2}\right)\notag\\
		& \left.+\sum_{l \in \mathbb{Z}}\left|\int_{0}^{\frac{1}{4 N}}\mathcal{F}(\mathbf{f})_{2,2}(x) e^{-2 \pi i(4 N) lx} d x\right|^{2}+\sum_{l \in \mathbb{Z}}\left|\int_{0}^{\frac{1}{4 N}}\mathcal{F}(\mathbf{f})_{2,2}(x) e^{-4 \pi i r x} e^{-2 \pi i(4 N) lx} dx\right|^{2}\right).
	\end{align}
	By applying Lemma \ref{RefH2}, we have
	\begin{align}\label{eqotto2}
		\sum_{\substack{ j=1 \notag\\
				\lambda \in \Lambda}}^2\left|\left\langle\mathbf{f},\mathcal{D}_{2 N \lambda} \mathbf{f}_{j}\right\rangle\right|^{2}= & \frac{1}{2}\left(\int_{0}^{\frac{1}{4 N}}\left|\mathcal{F}(\mathbf{f})_{1,1}(x)\right|^{2} dx+\int_{0}^{\frac{1}{4 N}}\left|\mathcal{F}(\mathbf{f})_{1,1}(x)\right|^{2} dx\right. \notag\\
		& +\int_{0}^{\frac{1}{4 N}}\left|\mathcal{F}(\mathbf{f})_{1,2}(x)\right|^{2} dx+\int_{0}^{\frac{1}{4 N}}\left|\mathcal{F}(\mathbf{f})_{1,2}(x)\right|^{2} dx \notag\\
		& +\int_{0}^{\frac{1}{4 N}}\left|\mathcal{F}(\mathbf{f})_{2,1}(x)\right|^{2} dx+\int_{0}^{\frac{1}{4 N}}\left|\mathcal{F}(\mathbf{f})_{2,1}(x)\right|^{2} dx \notag \\
		& \left.+\int_{0}^{\frac{1}{4 N}}\left|\mathcal{F}(\mathbf{f})_{2,2}(x)\right|^{2} dx+\int_{0}^{\frac{1}{4 N}}\left|\mathcal{F}(\mathbf{f})_{2,2}(x)\right|^{2} dx\right)\notag \\
		= & \int_{0}^{\frac{1}{4 N}}\left|\mathcal{F}(\mathbf{f})_{1,1}(x)\right|^{2} dx+\int_{0}^{\frac{1}{4 N}}\left|\mathcal{F}(\mathbf{f})_{1,2}(x)\right|^{2} dx \notag\\
		&+\int_{0}^{\frac{1}{4 N}}\left|\mathcal{F}(\mathbf{f})_{2,1}(x)\right|^{2} dx+\int_{0}^{\frac{1}{4 N}}\left|\mathcal{F}(\mathbf{f})_{2,2}(x)\right|^{2} dx.
	\end{align}
	Now, by selecting $\mathbf{f}_{t} \in \ell^{2}\left(\Lambda,\mathcal{M}_{2}\right)$ in accordance with its Fourier transformation,
	 \begin{align}\label{eat1}
	\mathcal{F}\left(\mathbf{f}_{t}\right)(x)=\left[\begin{array}{cc}
		\phi_{\left[0, \frac{1}{4 N}[\right.}+\frac{1}{a_{0}} \phi_{\left[\frac{1}{4 N}, \frac{2}{4 N}\right]}  &\phi_{\left[0, \frac{1}{4 N}[\right.}+\frac{1}{a_{0}} \phi_{\left[\frac{1}{4 N}, \frac{2}{4 N}\right]} \\
		\phi_{\left[0, \frac{1}{4 N}[\right.}+\frac{1}{a_{0}} \phi_{\left[\frac{1}{4 N}, \frac{2}{4 N}\right]} &\phi_{\left[0, \frac{1}{4 N}[\right.}+\frac{1}{a_{0}} \phi_{\left[\frac{1}{4 N}, \frac{2}{4 N}\right]}
	\end{array}\right].
	\end{align}
	Upon substituting $\mathcal{F}\left(\mathbf{f}_{t}\right)(x)$ into equation \eqref{eqotto2} ,we have
	\begin{align} \label{equat1}
	\sum_{\substack{ j=1 \\
	\lambda \in \Lambda}}^2\left|\left\langle\mathbf{f},\mathcal{D}_{2 N \lambda} \mathbf{f}_{j}\right\rangle\right|^{2}=\frac{1}{N}.
\end{align}
Then by using \eqref{euu12}, \eqref{eat1} and \eqref{equat1}, we ultimately arrive at a contradiction.
	\end{remark}	
The following lemma will be used in the characterization of matrix-valued DFNS.
	\begin{lemma}\label{Blem12}
For each $j \in \mathbb{N}_p$, let $\mathbf{0}\ne\mathbf{f}_{j} \in  \ell^{2}\left(\Lambda,\mathcal{M}_{n}\right)$. Suppose that
\begin{enumerate}
\item   $\left\|\mathcal{F}\left( \mathbf{f}_{j}\right)(x)\right\|_{\mathcal{M}_{n}} \leq  b_{0},  \, a.e. \ x \in \Omega$, $j \in \mathbb{N}_p$, where $b_{0}$  is a positive real number.\label{hypo12}
\item
	$\Upsilon_{j,m,k}(x) =
	\begin{bmatrix}
		\upsilon_1(x) \\
		\upsilon_2(x)
	\end{bmatrix}$,
	where
	\begin{align*}
	{\upsilon}_1(x) =
	\begin{bmatrix}
		\mathcal{F}(\mathbf{f}_j)_{m,k}(x) \\
		\mathcal{F}(\mathbf{f}_j)_{m,k}\left(x + \frac{1}{4N}\right) \\
		\vdots \\
		\mathcal{F}(\mathbf{f}_j)_{m,k}\left(x + \frac{2N-1}{4N}\right)
	\end{bmatrix}_{2N \times 1}
	\quad
	\text{and}
\quad
	{\upsilon}_2(x) =
	\begin{bmatrix}
	\mathcal{F}(\mathbf{f}_j)_{m,k}\left(x + \frac{N}{2}\right) \\
	\mathcal{F}(\mathbf{f}_j)_{m,k}\left(x + \frac{N}{2} + \frac{1}{4N}\right) \\
	\vdots \\
	\mathcal{F}(\mathbf{f}_j)_{m,k}\left(x + \frac{N}{2} + \frac{2N-1}{4N}\right)
	\end{bmatrix}_{2N \times 1}.
\end{align*}\label{hypo13}
\item 	
${\Psi}_{j,m,k}(x) = \mathbf{E}(x) \odot \Upsilon_{j,m,k}(x)$ = Hadamard product (entry-wise multiplication) of $\Upsilon_{j,m,k}(x)$ with
$\mathbf{E}(x) =
	\begin{bmatrix}
		\mathbf{E}_1(x) \\
		\mathbf{E}_1(x)
	\end{bmatrix}$, where
$\mathbf{E}_1(x) =
	\begin{bmatrix}
		e^{4 \pi i r x} \\
		e^{4 \pi i r \left(x + \frac{1}{4N}\right)} \\
		\vdots \\
		e^{4 \pi i r \left(x + \frac{2N-1}{4N}\right)}
	\end{bmatrix}_{2N \times 1}$\label{hypo14}.
\item for $ m,k \in\{1,2, \ldots, n\}$ and $x \in\left[0, \frac{1}{4 N}\right)$, define a $4 N \times 2p$ matrix as follows:
 \begin{align}\label{tam99}
\Gamma_{m, k}(x)=
 \left[\begin{array}{llllll}
{\Upsilon}_{1, m, k}(x) & {\Psi}_{1, m, k}(x) & {\Upsilon}_{2, m, k}(x) & {\Psi}_{2, m, k}(x) \cdots & {\Upsilon}_{p, m, k}(x) & \left.{\Psi}_{p, m, k}(x)\right]_{4 N \times 2p}.
\end{array}\right.
\end{align}
\end{enumerate}
Then, for any $\mathbf{f} \in \ell^{2}\left(\Lambda,\mathcal{M}_{n}\right)$, we have
\begin{align*}
	4 N \left\| \left\{ \left\langle \mathbf{f}, \mathcal{D}_{2N\lambda} \mathbf{f}_{j} \right\rangle \right\}_{\substack{\lambda \in \Lambda \\  j \in \mathbb{N}_p,}} \right\|_{\ell^2}^2=\int_{0}^{\frac{1}{4 N}}\left\|\sum_{m,k=1}^{n} \Gamma_{m, k}^{*}(x) \mathscr{G}_{\mathbf{f}_{m,k}}(x)\right\|^{2} dx,
\end{align*}
	where $\Gamma_{m, k}^{*}(x)$ is conjugate transpose of matrix $\Gamma_{m, k}(x)$
	and
	$\mathscr{G}_{\mathbf{f}_{m,k}}(x) =
	\begin{bmatrix}
		\mathscr{G}_1(x) \\
		\mathscr{G}_2(x)
	\end{bmatrix}
	$,
	where
	\begin{align*}
		\mathscr{G}_1(x) =
		\begin{bmatrix}
			\mathcal{F}(\mathbf{f})_{m,k}(x)  \\
			\mathcal{F}(\mathbf{f})_{m,k}\left(x+\frac{1}{4 N}\right) \\
			\vdots \\
			\mathcal{F}(\mathbf{f})_{m,k}\left(x+\frac{2 N-1}{4 N}\right)
		\end{bmatrix}_{2N \times 1}
		\quad \text{and} \quad
		\mathscr{G}_2(x) =
		\begin{bmatrix}
			\mathcal{F}(\mathbf{f})_{m,k}\left(x+\frac{N}{2}\right) \\
			\mathcal{F}(\mathbf{f})_{m,k}\left(x+\frac{N}{2}+\frac{1}{4 N}\right) \\
			\vdots \\
			\mathcal{F}(\mathbf{f})_{m,k}\left(x+\frac{N}{2}+\frac{2 N-1}{4 N}\right)
		\end{bmatrix}_{2N \times 1}.
	\end{align*}
\end{lemma}
\proof
	Consider  $\mathbf{f} \in \ell^{2}\left(\Lambda, \mathcal{M}_{n}\right)$ be arbitrary. For $j \in \mathbb{N}_p$ and $x \in\left[0, \frac{1}{2}\right)$, write
	\begin{align}\label{Eqq1}
		\mathfrak{X}_j(x)
		=  \sum_{m,k=1}^{n}\left(\mathcal{F}(\mathbf{f})_{m,k}(x) \overline{\mathcal{F}\left( \mathbf{f}_{j}\right)_{m,k}(x)}\right.		\left.+\mathcal{F}(\mathbf{f})_{m,k}\left(x+\frac{N}{2}\right)\overline{\mathcal{F}\left(\mathbf{f}_{j}\right)_{m,k}\left(x+\frac{N}{2}\right)}\right).
	\end{align}
Now, for $x \in\left[0, \frac{1}{4 N}\right)$, we  have
	\begin{align*}
		\sum_{m,k=1}^{n} \Gamma_{m,k}^*(x) \, \mathscr{G}_{\mathbf{f}_{m,k}}(x)
		=
		\begin{bmatrix}
			\sum\limits_{m,k=1}^{n} \left\langle \Upsilon_{1,m,k}(x), \mathscr{G}_{\mathbf{f}_{m,k}}(x) \right\rangle \\
			\sum\limits_{m,k=1}^{n} \left\langle \Psi_{1,m,k}(x), \mathscr{G}_{\mathbf{f}_{m,k}}(x) \right\rangle \\
			\vdots \\
			\sum\limits_{m,k=1}^{n} \left\langle \Upsilon_{p,m,k}(x), \mathscr{G}_{\mathbf{f}_{m,k}}(x) \right\rangle \\
			\sum\limits_{m,k=1}^{n} \left\langle \Psi_{p,m,k}(x), \mathscr{G}_{\mathbf{f}_{m,k}}(x) \right\rangle
		\end{bmatrix}_{2p \times 1}.
	\end{align*}
Using hypotheses \ref{hypo13} and \ref{hypo14}, the inner products in above matrix are equal to
	\begin{align*}
		\left\langle \Upsilon_{j,m,k}(x), \mathscr{G}_{\mathbf{f}_{m,k}}(x) \right\rangle &=
		\sum_{g=0}^{2N-1}
		\overline{\mathcal{F}(\mathbf{f}_j)_{m,k}\left(x + \frac{g}{4N}\right)} \cdot
		\mathcal{F}(\mathbf{f})_{m,k}\left(x + \frac{g}{4N}\right)\\
		& +
		\sum_{g=0}^{2N-1}
		\overline{\mathcal{F}(\mathbf{f}_j)_{m,k}\left(x + \frac{N}{2} + \frac{g}{4N}\right)} \cdot
		\mathcal{F}(\mathbf{f})_{m,k}\left(x + \frac{N}{2} + \frac{g}{4N}\right),
	\end{align*}
and
	\begin{align*}
		&\left\langle \Psi_{j,m,k}(x), \mathscr{G}_{\mathbf{f}_{m,k}}(x) \right\rangle\\
		&= \sum_{g=0}^{2N - 1}
		\overline{e^{4 \pi i r \left(x + \frac{g}{4N} \right)} \, \mathcal{F}(\mathbf{f}_j)_{m,k} \left(x + \frac{g}{4N} \right)}
		\cdot \mathcal{F}(\mathbf{f})_{m,k} \left(x + \frac{g}{4N} \right) \\
		&\quad + \sum_{g=0}^{2N - 1}
		\overline{e^{4 \pi i r \left(x + \frac{N}{2} + \frac{g}{4N} \right)} \, \mathcal{F}(\mathbf{f}_j)_{m,k} \left(x + \frac{N}{2} + \frac{g}{4N} \right)} \cdot \mathcal{F}(\mathbf{f})_{m,k} \left(x + \frac{N}{2} + \frac{g}{4N} \right).
	\end{align*}
Therefore, using equation \eqref{Eqq1}, we get
\begin{align*}
	\begin{gathered}
		\sum_{m,k=1}^{n} \Gamma_{m, k}^{*}(x) \mathscr{G}_{\mathbf{f}_{m,k}}(x)=\left[\begin{array}{c}
			\sum_{g=0}^{2 N-1} \mathfrak{X}_1
			\left(x+\frac{g}{4 N}\right)\\
			\sum_{g=0}^{2 N-1} e^{-4 \pi i r\left(x+\frac{g}{4 N}\right)}  \mathfrak{X}_1
			\left(x+\frac{g}{4 N}\right)\\
			\sum_{g=0}^{2 N-1} \mathfrak{X}_2
			\left(x+\frac{g}{4 N}\right)\\
			\sum_{g=0}^{2 N-1}e^{-4 \pi i r\left(x+\frac{g}{4 N}\right)} \mathfrak{X}_2\left(x+\frac{g}{4 N}\right) \\
			\vdots\\
			\sum_{g=0}^{2 N-1} \mathfrak{X}_p
			\left(x+\frac{g}{4 N}\right)\\
			\sum_{g=0}^{2 N-1}e^{-4 \pi i r\left(x+\frac{g}{4 N}\right)} \mathfrak{X}_p\left(x+\frac{g}{4 N}\right)\\
		\end{array}\right]_{2p \times 1},
	\end{gathered}
\end{align*}
and
\begin{align}\label{eqint1}
&\int\limits_{0}^{\frac{1}{4N}}\left\|\sum_{m,k=1}^{n} \Gamma_{m, k}^{*}(x) \mathscr{G}_{\mathbf{f}_{m,k}}(x)\right\|_{\mathbb{C}^{2p}}^{2} \, dx \notag\\
		 &=\sum_{j=1 }^{p}\int\limits_{0}^{\frac{1}{4N}}\left(\left|\sum_{g=0}^{2N-1} \mathfrak{X}_j\left(x+\frac{g}{4N}\right)\right|^{2}+\left|\sum_{g=0}^{2N-1} e^{-\pi i r \frac{g}{N}} \mathfrak{X}_j\left(x+\frac{g}{4N}\right)\right|^{2}\right) \, dx.
\end{align}
For any $\mathbf{f} \in \ell^{2}\left(\Lambda, \mathcal{M}_{n}\right)$, we compute
	\begin{align} \label{euu1}
		\left\| \left\{ \left\langle \mathbf{f}, \mathcal{D}_{2N\lambda} \mathbf{f}_{j} \right\rangle \right\}_{\substack{\lambda \in \Lambda \\  j \in \mathbb{N}_p}} \right\|_{\ell^2}^2
& = \sum_{\substack{j=1 \atop \lambda \in \Lambda}}^{p}\left|\left\langle\mathbf{f}, \mathcal{D}_{2 N \lambda} \mathbf{f}_{j}\right\rangle\right|^{2} \notag \\
		= & \sum_{\substack{ j=1\\
				\lambda \in \Lambda}}^{p}\left|\left\langle\mathcal{F}(\mathbf{f})(x), e^{4 \pi i N \lambda\left(x\right)} \mathcal{F}\left( \mathbf{f}_{j}\right)(x)\right\rangle\right|^{2} \notag\\
		= & \sum_{\substack{j=1\\
				\lambda \in \Lambda}}^{p}\left|\int_{\Omega} \sum_{m,k=1}^{n}\mathcal{F}(\mathbf{f})_{m,k}(x) \overline{\mathcal{F}\left( \mathbf{f}_{j}\right)_{m,k}(x)} e^{-4 \pi i N \lambda\left(x\right)} d x\right|^{2} \notag \\
		= & \sum_{\substack{ j=1 \\
				l \in \mathbb{Z}}}^{p}\left|\int_{\Omega} \sum_{m,k=1}^{n}\mathcal{F}(\mathbf{f})_{m,k}(x)  \overline{\mathcal{F}\left(\mathbf{f}_{j}\right)_{m,k}(x)} e^{-4 \pi i N(2 l)\left(x\right)} d x\right|^{2} \notag \\
		+ & \sum_{\substack{ j=1 \\
				l \in \mathbb{Z}}}^{p}\left|\int_{\Omega} \sum_{m,k=1}^{n}\mathcal{F}(\mathbf{f})_{m,k}(x) \overline{\mathcal{F}\left(\mathbf{f}_{j}\right)_{m,k}(x)} e^{-4 \pi i N\left(\frac{r}{N}+2 l\right)\left(x\right)} d x\right|^{2} .
	\end{align}
	Using \eqref{Eqq1} and  \eqref{euu1}, we have
	\begin{align}\label{eqq123}
		 \sum_{\substack{j=1\\
				\lambda \in \Lambda}}^{p}\left|\left\langle\mathbf{f}, \mathcal{D}_{2 N \lambda} \mathbf{f}_{j}\right\rangle\right|^{2}
		& =\sum_{\substack{ j=1\\
				l \in \mathbb{Z}}}^{p}\left|\int_{0}^{\frac{1}{2}} \mathfrak{X}_j(x) e^{-4 \pi i N(2 l)\left(x\right)} d x\right|^{2} \notag \\
		& +\sum_{\substack{ j=1\\
				l \in \mathbb{Z}}}^{p}\left|\int_{0}^{\frac{1}{2}} \mathfrak{X}_j(x)
		e^{-4 \pi i N\left(\frac{r}{N}+2 l\right)\left(x\right)} d x\right|^{2} \notag \\
		& =\sum_{\substack{ j=1\\
				l \in \mathbb{Z}}}^{p}\left|\int_{0}^{\frac{1}{4 N}} \sum_{g=0}^{2 N-1} \mathfrak{X}_j\left(x+\frac{g}{4 N}\right) e^{-2 \pi i(4 N) l\left(x\right)} d x\right|^{2} \notag \\
		& +\sum_{\substack{ j=1\\
				l \in \mathbb{Z}}}^{p}\left|\int_{0}^{\frac{1}{4 N}} \sum_{g=0}^{2 N-1} \mathfrak{X}_j\left(x+\frac{g}{4 N}\right) e^{-4 \pi i r\left(x+\frac{g}{4 N}\right)} e^{-2 \pi i(4 N) l\left(x\right)} d x\right|^{2} .
	\end{align}
	By using  Lemma \ref{RefH2},  hypothesis \ref{hypo12}  and \eqref{eqq123}, we have
\begin{align}\label{eu9}
&\left\| \left\{ \left\langle \mathbf{f}, \mathcal{D}_{2N\lambda} \mathbf{f}_{j} \right\rangle \right\}_{\substack{\lambda \in \Lambda \\  j \in \mathbb{N}_p}} \right\|_{\ell^2}^2 \notag \\
&=\frac{1}{4N} \Bigg(\sum_{j=1}^{p} \int\limits_{0}^{\frac{1}{4N}} \left| \sum_{g=0}^{2N-1}\mathfrak{X}_j\left(x + \frac{g}{4N}\right) \right|^2 dx
+ \sum_{j=1}^{p} \int\limits_{0}^{\frac{1}{4N}} \left| \sum_{g=0}^{2N-1} \mathfrak{X}_j\left(x + \frac{g}{4N}\right) e^{-i \pi r \frac{g}{N} } \right|^2 dx \Bigg).
\end{align}
The desired result follows from equations \eqref{eqint1} and \eqref{eu9} .
\endproof
Now, we are ready for necessary and sufficient conditions for the existence of matrix-valued DFNS.
\begin{theorem}\label{MainthmI}
 Suppose
 \begin{enumerate}
\item $\big\{\mathbf{f}_j\big\}_{j \in \mathbb{N}_p} \subset \ell^{2}\left(\Lambda,\mathcal{M}_{n}\right)$ satisfies
	$\left\| \mathcal{F} \left( \mathbf{f}_j \right)(x) \right\|_{\mathcal{M}_n} \leq b_0$ for each  $j \in \mathbb{N}_p$,  a.e. $x \in \Omega$,  where $b_{0}$ is a positive real number.\label{Hypot1}
\item ${\Gamma}_{m, k}(x)$ is a  $4N \times 2p$ matrix as defined in \eqref{tam99} for  $ m,k \in\{1,2, \ldots, n\}$, and  a.e.  $x \in\left[0, \frac{1}{4 N}\right)$.
\end{enumerate}
Then,  $ \left\{ \mathcal{D}_{2 N \lambda} \mathbf{f}_{j}\right\}_{\lambda\in \Lambda \atop j \in \mathbb{N}_p}$  is a matrix-valued DFNS for  $\ell^{2}\left(\Lambda,\mathcal{M}_{n}\right)$ with  lower frame bound $a_{0}$ if and only if
	\begin{align}\label{Hy324}
		\left\|\sum_{m,k=1}^{n} {\Gamma}_{m, k}^{*}(x) \mathbf{w}_{m,k}\right\|_{\mathbb{C}^{2p}}^{2} \geq 4 N a_{0}\sum_{m,k=1}^{n}\left\|\mathbf{w}_{m,k}\right\|_{\mathbb{C}^{4 N}}^{2} \text {, a.e. } x \in\left(0, \frac{1}{4 N}\right),
\end{align}
	for all $\mathbf{w}_{m,k} \in \mathbb{C}^{4 N}$,  $1\leq m, k\leq  n$.
\end{theorem}
\proof
 Given that $\left\|\mathcal{F}\left( \mathbf{f}_{j}\right)(x)\right\|_{\mathcal{M}_{n}} \leq b_{0}<\infty$, for $ j \in \mathbb{N}_p$, a.e. $ x\in \Omega$, so by Lemma \ref{Blem12}, for each $\mathbf{f} \in \ell^{2}\left(\Lambda,\mathcal{M}_{n}\right)$, we have
		\begin{align}\label{Bre3250}
		4 N \left\| \left\{ \left\langle \mathbf{f}, \mathcal{D}_{2N\lambda} \mathbf{f}_{j} \right\rangle \right\}_{\substack{\lambda \in \Lambda \\ j \in \mathbb{N}_p}} \right\|_{\ell^2}^2=\int_{0}^{\frac{1}{4 N}}\left\|\sum_{m,k=1}^{n} \Gamma_{m, k}^{*}(x) \mathscr{G}_{\mathbf{f}_{m,k}}(x)\right\|^{2} dx,
	\end{align}
	where \( \mathscr{G}_{\mathbf{f}_{m,k}}(x) \) is defined as in Lemma \ref{Blem12}.
 Assume  that \( \left\{ \mathcal{D}_{2 N \lambda} \mathbf{f}_{j}\right\}_{\lambda\in \Lambda \atop j \in \mathbb{N}_p} \)
	is a matrix-valued DFNS for \( \ell^{2}\left(\Lambda,\mathcal{M}_{n}\right) \) with lower frame bound \( a_{0} \)(say).
	Let \( f \in L^{2}\left(0, \frac{1}{4 N}\right) \) be an arbitrary element.
	For \( m,k \in\{1,2, \ldots, n\} \), let
\begin{align}\label{Bre326}
	\mathbf{w}_{m,k}=
	\begin{bmatrix}
		({c}_{m,k})_{1} \\
		({c}_{m,k})_{2} \\
		\vdots \\
		({c}_{m,k})_{2N} \\
		({c}_{m,k})_{2N+1} \\
		\vdots \\
		({c}_{m,k})_{4N}
	\end{bmatrix}_{4N \times 1} \in \mathbb{C}^{4N}.
\end{align}
	Describe $\mathbf{f} \in \ell^{2}\left(\Lambda, \mathcal{M}_{n}\right)$
	in terms of its Fourier transform as:
	\[
	\mathcal{F}(\mathbf{f}(x)) =
	\begin{bmatrix}
		\mathcal{F}(\mathbf{f})_{1,1}(x) & \mathcal{F}(\mathbf{f})_{1,2}(x) & \cdots & \mathcal{F}(\mathbf{f})_{1,n}(x) \\
		\mathcal{F}(\mathbf{f})_{2,1}(x)& \mathcal{F}(\mathbf{f})_{2,2}(x) & \cdots & \mathcal{F}(\mathbf{f})_{2,n}(x)\\
		\vdots & \vdots & \ddots & \vdots \\
		\mathcal{F}(\mathbf{f})_{n,1}(x) & \mathcal{F}(\mathbf{f})_{n,2}(x)) & \cdots & \mathcal{F}(\mathbf{f})_{n,n}(x)
	\end{bmatrix}
	\in L^{2}\left(\Omega,\mathcal{M}_{n}\right).
	\]
	For any fixed $m_ {0},k_{0} \in\{1,2, \ldots, n\}$, we have following cases:\\
	$\mathbf{Case 1: }$ If $x \in\left[0, \frac{1}{2}\right)$, then there exists unique $t^{\prime} \in\{0,1, \ldots, 2 N-$ $1\}$ such that $x-\frac{t^{\prime}}{4 N} \in\left[0, \frac{1}{4 N}\right)$, so define
	$$
	\mathcal{F}(\mathbf{f})_{m_{0},k_{0}}(x)=({c}_{m_{0},k_{0}})_{t^{\prime}+1} f\left(x-\frac{t^{\prime}}{4 N}\right).
	$$
	$\mathbf{Case 2:}$ If $ x\in\left[\frac{N}{2}, \frac{N+1}{2}\right)$, there exists $t^{\prime \prime} \in\{0,1, \ldots, 2 N-1\}$ such that $x-\frac{t^{\prime \prime}}{4 N}-\frac{N}{2} \in\left[0, \frac{1}{4 N}\right)$, so define
	$$
	\mathcal{F}(\mathbf{f})_{m_{0},k_{0}}(x)=({c}_{m_{0},k_{0}})_{2 N+1+t^{\prime \prime}} f\left(x-\frac{t^{\prime \prime}}{4 N}-\frac{N}{2}\right).
	$$
	Thus, for $m,k \in\{1,2, \ldots, n\}, t \in\{0,1, \ldots, 2 N-1\}$ and $x \in\left[0, \frac{1}{4 N})\right.$, we have
	\begin{align}
		{\mathcal{F}(\mathbf{f})_{m,k}\left(x+\frac{t}{4 N}\right) } & =({c}_{m,k})_{t+1} f(x),  \label{Bre327}
\intertext{and}
		{\mathcal{F}(\mathbf{f})_{m,k}\left(x+\frac{t}{4 N}+\frac{N}{2}\right) } & =({c}_{m,k})_{2 N+t+1} f(x) \label{Bre328}.
	\end{align}
	Now, for $m,k \in\{1,2, \ldots, n\}$ and $x \in\left[0, \frac{1}{4 N}\right)$, we define  $\mathscr{G}_{\mathbf{f}_{m,k}}(x) =
	\begin{bmatrix}
		\mathscr{G}_1(x) \\
		\mathscr{G}_2(x)
	\end{bmatrix} \label{gamma1}
	$,
	where
	\begin{align}\label{Bre329}
		\mathscr{G}_1(x) =
		\begin{bmatrix}
			\mathcal{F}(\mathbf{f})_{m,k}(x)  \\
			\mathcal{F}(\mathbf{f})_{m,k}\left(x+\frac{1}{4 N}\right) \\
			\vdots \\
			\mathcal{F}(\mathbf{f})_{m,k}\left(x+\frac{2 N-1}{4 N}\right)
		\end{bmatrix}_{2N \times 1}
		\quad \text{and} \quad
		\mathscr{G}_2(x) =
		\begin{bmatrix}
			\mathcal{F}(\mathbf{f})_{m,k}\left(x+\frac{N}{2}\right) \\
			\mathcal{F}(\mathbf{f})_{m,k}\left(x+\frac{N}{2}+\frac{1}{4 N}\right) \\
			\vdots \\
			\mathcal{F}(\mathbf{f})_{m,k}\left(x+\frac{N}{2}+\frac{2 N-1}{4 N}\right)
		\end{bmatrix}_{2N \times 1}.
	\end{align}
	By \eqref{Bre326}, \eqref{Bre327}, \eqref{Bre328} and \eqref{Bre329}, for $m,k \in\{1,2, \ldots, n\}$ and $ x\in\left[0, \frac{1}{4 N}\right)$, we have
	\begin{align}\label{Bre330}
		\mathscr{G}_{\mathbf{f}_{m,k}}(\xi)=\mathbf{w}_{m,k} f(x).
	\end{align}
	Now, using \eqref{Bre330} and \eqref{Bre3250}, we have
	\begin{align}\label{Bre331}
	\int_{0}^{\frac{1}{4 N}}|f(x)|^{2} \left\|\sum_{m,k=1}^{n} {\Gamma}_{m, k}^{*}(x) \mathbf{w}_{m,k}\right\|_{\mathbb{C}^{2p}} d x & =\int_{0}^{\frac{1}{4 N}}\left\|\sum_{m,k=1}^{n} {\Gamma}_{m, k}^{*}(x) \mathscr{G}_{\mathbf{f}_{m,k}}(x)\right\|_{\mathbb{C}^{2P}}^{2} d x \notag\\
	& =4 N \left\| \left\{ \left\langle \mathbf{f}, \mathcal{D}_{2N\lambda} \mathbf{f}_{j} \right\rangle \right\}_{\substack{\lambda \in \Lambda \\ j \in \mathbb{N}_p}} \right\|_{\ell^2}^2.
	\end{align}
	Therefore, by \eqref{Bre331} and the hypothesis that $ \left\{ \mathcal{D}_{2 N \lambda} \mathbf{f}_{j}\right\}_{\lambda\in \Lambda \atop j \in \mathbb{N}_p}$  is a frame with lower frame bound $a_{0}$, we get
\begin{align}\label{Bre332}
& \int_{0}^{\frac{1}{4 N}}|f(x)|^{2} \left\|\sum_{m,k=1}^{n} {\Gamma}_{m, k}^{*}(x) \mathbf{w}_{m,k}\right\|_{\mathbb{C}^{2p}}^{2} d x \notag\\
& \geq 4 N a_{0}\|\mathbf{f}\|_{\ell^{2}}^{2} \notag\\
& =4 N a_{0}\|\mathcal{F}(\mathbf{f})\|_{L^{2}}^{2} \notag\\
& =4 N a_{0} \sum_{m,k=1}^{n} \sum_{g=0}^{2 N-1}\left(\int\limits_{0}^{4 N}\left(\left|\mathcal{F}(\mathbf{f})_{m,k}\left(x+\frac{g}{4 N}\right)\right|^{2}+\left|\mathcal{F}(\mathbf{f})_{m,k}\left(x+\frac{g}{4 N}+\frac{N}{2}\right)\right|^{2}\right) d x.\right.
\end{align}
	Using \eqref{Bre327}, \eqref{Bre328} and \eqref{Bre332}, we have
\begin{align*}
	\int_{0}^{\frac{1}{4 N}}|f(x)|^{2}\left\|\sum_{m,k=1}^{n} {\Gamma}_{m, k}^{*}(x) \mathbf{w}_{m,k}\right\|_{\mathbb{C}^{2p}}^{2} d x \geq 4 N a_{0} \sum_{m,k=1}^{n} \int_{0}^{\frac{1}{4 N}}\left\|\mathbf{w}_{m,k}\right\|_{\mathbb{C}^{4 N}}^{2}|f(x)|^{2} dx.
\end{align*}
	This gives
\begin{align*}
	\left\|\sum_{m,k=1}^{n} {\Gamma}_{m, k}^{*}(x) \mathbf{w}_{m,k}\right\|_{\mathbb{C}^{2p}}^{2} \geq 4 N a_{0} \sum_{m,k=1}^{n}\left\|\mathbf{w}_{m, k}\right\|_{\mathbb{C}^{4 N}}^{2} \  \text{a.e.} \  x \in\left(0, \frac{1}{4 N}\right).
\end{align*}
Thus, \eqref{Hy324} is proved.

	Conversely , assume that \eqref{Hy324}  holds. Using condition \eqref{Hypot1},
	$\left\|\mathcal{F}\left(\mathbf{f}_{j}\right)(x)\right\|_{\mathcal{M}_{n}} \leq b_{0}$,  for  $j \in \mathbb{N}_p$,  a.e. $x \in \Omega$,
	so by Theorem \ref{thmM1}, the collection
	$\left\{ \mathcal{D}_{2 N \lambda} \mathbf{f}_{j} \right\} _{\lambda\in \Lambda \atop j \in \mathbb{N}_p}
	$
	forms a  Bessel sequence with Bessel bound $2^{p-1} b_{0}^2n^2$. Hence, we only need to verify the lower frame inequality. From \eqref{Bre3250}, for each
	$
	\mathbf{f} \in \ell^{2}\left(\Lambda ,\mathcal{M}_{n}\right),
	$
	we have
	\begin{align}\label{equa123}
	4 N \sum_{\substack{ j=1 \\ \lambda \in \Lambda}}^{p}\left|\left\langle\mathbf{f}, \mathcal{D}_{2 N \lambda} \mathbf{f}_{j}\right\rangle\right|^{2}=\int_{0}^{\frac{1}{4 N}}\left\|\sum_{m,k=1}^{n} \Gamma_{m, k}^{*}(x) \mathscr{G}_{\mathbf{f}_{m,k}}(x)\right\| dx,
	\end{align}
	where $\mathscr{G}_{\mathbf{f}_{m,k}}(x)$  as defined in \eqref{gamma1}. Therefore, by \eqref{Hy324} and \eqref{equa123}, we have
	\begin{align*}
		\sum_{\substack {j=1 \\ \lambda \in \Lambda}}^{p} \left|\left\langle \mathbf{f}, \mathcal{D}_{2 N \lambda} \mathbf{f}_{j} \right\rangle \right|^{2}
		 \geq a_{0} \int_{0}^{\frac{1}{4 N}} \sum_{m,k=1}^{n} \left\|\mathscr{G}_{\mathbf{f}_{m,k}}(x) \right\|_{\mathbb{C}^{4 N}}^{2} \, d\xi
		 = a_{0} \|\mathcal{F}(\mathbf{f})\|_{L^{2}}^{2}
		 = a_{0} \|\mathbf{f}\|_{\ell^{2}}^{2}.
\end{align*}
This completes the proof.
\endproof
We conclude the section with an applicative example of Theorem \ref{MainthmI}.
\begin{example}\label{Ex3.10}
 Let $N=2$, $r=1$, $n=2$ and $p=8$. Then, $\Lambda = \left\{ 0, \frac{1}{2} \right\} + 2 \mathbb{Z}$ and
	$\Omega = \left[ 0, \frac{1}{2} \right) \cup \left[ 1, \frac{3}{2} \right)$
	and  consider  $\{\mathbf{f_j} \}_{j \in \mathbb{N}_{8}}  \subset \ell^{2}\left(\Lambda, \mathcal{M}_{2}\right)$ as in Example \ref{Exam1}. Then,
\begin{align*}
	\left\|\mathcal{F}\left( \mathbf{f}_{j}\right)(x)\right\|_{\mathcal{M}_{n} }\leq 2, \text { a.e. } x \in \Omega, \, j \in \mathbb{N}_{8},
\end{align*}
 and the collection $ \left\{ \mathcal{D}_{2 N \lambda} \mathbf{f}_{j}\right\}_{\lambda\in \Lambda \atop j \in \mathbb{N}_{8}}$  is a Bessel sequence with  Bessel bound $ 2^{11}$.
		For $m,k  \in\{1,2\}$, let $\Gamma_{m, k}(x)$ be the $8 \times 16$ matrices as defined in \eqref{tam99} . By utilizing  any matrix calculation software, it is straightforward to verify that for each $m,k \in \{1,2\}$, we have
		\begin{align}\label{exavv}
		{\Gamma}_{m, k}(x){\Gamma}_{m^{\prime}, k^{\prime}}^{*}(x)=8 \delta_{k k^{\prime}}\delta_{m m^{\prime}} \mathbf{I}_{8}  \text { for } 1 \leq m,m^{\prime },k, k^{\prime} \leq 2.
	\end{align}
	For each $\mathbf{w}_{1,1}, \mathbf{w}_{1,2},\mathbf{w}_{2,1},\mathbf{w}_{2,2}\in \mathbb{C}^{8}$ and a.e. $x \in\left(0, \frac{1}{8}\right)$, using \eqref{exavv}, we compute
		\begin{align}\label{equ00}
		 \left\|\sum_{\substack{m,k=1}}^{2}  {\Gamma}_{m, k}^{*}(x) \mathbf{w}_{m,k}\right\|_{\mathbb{C}^{16}}^{2}
		& =\left\langle\sum_{\substack{m,k=1}}^{2}{\Gamma}_{m, k}^{*}(x) \mathbf{w}_{m,k}, \sum_{\substack{m^{\prime},k^{\prime}=1}}^{2} {\Gamma}_{m^{\prime}, k^{\prime}}^{*}(x) \mathbf{w}_{m^{\prime},k^{\prime}}\right\rangle \notag \\
		& =\left\langle 8 \mathbf{I}_{8} \mathbf{w}_{1,1}, \mathbf{w}_{1,1}\right\rangle+\left\langle 8 \mathbf{I}_{8} \mathbf{w}_{1,2}, \mathbf{w}_{1,2}\right\rangle+\left\langle 8 \mathbf{I}_{8} \mathbf{w}_{2,1}, \mathbf{w}_{2,1}\right\rangle+\left\langle 8 \mathbf{I}_{8} \mathbf{w}_{2,2}, \mathbf{w}_{2,2}\right\rangle \notag \\
		& =8\left(\left\|\mathbf{w}_{1,1}\right\|^{2}+\left\|\mathbf{w}_{1,2}\right\|^{2}+\left\|\mathbf{w}_{2,1}\right\|^{2}+\left\|\mathbf{w}_{2,2}\right\|^{2}\right) \notag \\
		& =8 \sum_{\substack{m,k=1}}^{2}\left\|\mathbf{w}_{m,k}\right\|^{2}.
	\end{align}
	Hence, by \eqref{equ00} and Theorem \ref{MainthmI},  $ \left\{ \mathcal{D}_{2 N \lambda} \mathbf{f}_{j}\right\}_{\lambda\in \Lambda \atop j \in \mathbb{N}_{8}}$ is a  matrix-valued DFNS for  $\ell^{2}\left(\Lambda, \mathcal{M}_{2}\right)$ with bounds $1$ and $2^{11}$.
 \end{example}
\section{Perturbation Results }
In signal processing, time-frequency analysis, quantum physics and duality relations etc.,
 it is important that fundamental properties of frames are kept under perturbation \cite{GosoI, BHan2, Heil2, KPSA, NPS}. Perturbation of frames and Riesz bases  were studied by Favier and Zalik in \cite{Favier}.  Most recently, by using a contraction between the affine group and the Weyl-Heisenberg group, perturbation of frames with Gabor and wavelet structure  studied by authors of \cite{JindIII}. Perturbation of discrete frames and their applications  in distributed signal processing can be found in \cite{deep2}. The deformation problem for small times is related to perturbation theory \cite{GosoI}.
 Jyoti and Vashisht  discussed perturbation  for dual frames of wave packets in the paper \cite{Jlk4}. Heil \cite{Heil1} is a good reference for perturbation of frames and Schauder bases in separable Banach spaces.
  In this section, we give perturbation of matrix-valued discrete frames of non-uniform shifts with explicit frame bounds. We begin with the perturbation of
matrix-valued DFNS in terms of Fourier transforms.
	\begin{theorem}\label{PerMI}
Suppose
\begin{enumerate}
  \item $\left\{ \mathcal{D}_{2 N \lambda} \mathbf{f}_{j}\right\}_{\lambda\in \Lambda \atop j \in \mathbb{N}_p}$  is matrix-valued DFNS for $\ell^{2}\left(\Lambda, \mathcal{M}_n\right)$ with frames bounds $a_{0}$ and $b_{0}$.
  \item  $\left\{\mathbf{g}_{j}\right\}_{j \in \mathbb{N}_{p}} \subset \ell^{2}\left(\Lambda,\mathcal{M}_{n}\right)$ satisfies
	\begin{align}
		\left\|\mathcal{F}
		\left(\mathbf{f}_{j}+\mathbf{g}_{j}\right)(x)\right\|_{\mathcal{M}_n }\leq \epsilon<2^{p-1} \epsilon^{2} n^{2}< a_{0}, \, j \in \mathbb{N}_p, \,  x\in \Omega. \label{Per4.1}
	\end{align}
\end{enumerate}
Then, the collection $\left\{ \mathcal{D}_{2 N \lambda} \mathbf{g}_{j} \right\}_{\lambda\in \Lambda \atop j \in \mathbb{N}_p}
	$ is also  a  matrix-valued DFNS for $\ell^{2}\left(\Lambda ,\mathcal{M}_n\right)$ with frame bounds $\left(\sqrt{a_{0}}-\sqrt{2^{p-1} \epsilon^{2} n^{2}}\right)^{2}$ and $\left(2^{p} \epsilon^{2} n^{2}+2b_{0}\right)$.
	\end{theorem}
	\proof Using  Lemma \ref{Blem1}, for each  $\mathbf{f} \in \ell^{2}\left(\Lambda,\mathcal{M}_{n}\right)$, we have
	\begin{align}
\sum_{\substack{ j=1 \\ \lambda \in \Lambda}}^{p}\left|\left\langle\mathbf{f}, \mathcal{D}_{2 N \lambda} \mathbf{g}_{j}\right\rangle\right|^{2}
		= & \sum_{\substack{ j=1\\
				\lambda \in \Lambda}}^{p}\left|\left\langle\mathbf{f}, \mathcal{D}_{2 N \lambda}\left(\mathbf{f}_{j}+\mathbf{g}_{j}\right)\right\rangle-\left\langle\mathbf{f},\mathcal{D}_{2 N \lambda} \mathbf{f}_{j}\right\rangle\right|^{2} \label{per4.3}\\
		\leq & 2\left(\sum_{\substack{ j=1 \\
				\lambda \in \Lambda}}^{p}\left|\left\langle\mathbf{f}, \mathcal{D}_{2 N \lambda}\left(\mathbf{f}_{j}+\mathbf{g}_{j}\right)\right\rangle\right|^{2}+\sum_{\substack{j=1 \\ \lambda \in \Lambda}}^{p}\left|\left\langle\mathbf{f},\mathcal{D}_{2 N \lambda} \mathbf{f}_{j}\right\rangle\right|^{2}\right) \label{per4.3zx}.
	\end{align}
From condition \eqref{Per4.1} and Theorem \ref{thmM1}, we have
	\begin{align}\label{Per4.4}
		\sum_{\substack{j=1 \\  \lambda \in \Lambda}}^{p}\left|\left\langle\mathbf{f},  \mathcal{D}_{2 N \lambda}\left(\mathbf{f}_{j}+\mathbf{g}_{j}\right)\right\rangle\right|^{2} \leq 2^{p-1} \epsilon^{2}n^{2}\|\mathbf{f}\|^{2},\ \mathbf{f} \in \ell^{2}\left(\Lambda,\mathcal{M}_{n}\right).
	\end{align}
From \eqref{per4.3zx} and \eqref{Per4.4} we have
\begin{align}\label{Per4.5} \sum\limits_{\substack{j=1 \\ \lambda \in \Lambda}}^{p}\left|\left\langle\mathbf{f}, \mathcal{D}_{2 N \lambda} \mathbf{g}_{j}\right\rangle\right|^{2} \leq\left(2^{p}\epsilon^{2}n^{2}+2b_{0}\right)\|\mathbf{f}\|^{2}.
\end{align}
Hence, the upper condition is obtained.
From  \eqref{per4.3} and \eqref{Per4.4}, we have
	\begin{align}\label{Per4.6}
	\sqrt{\sum_{\substack{j=1 \\ \lambda \in \Lambda}}^{p} \left| \left\langle \mathbf{f}, \mathcal{D}_{2 N \lambda} \mathbf{g}_{j} \right\rangle \right|^{2}}
		= & \sqrt{\sum_{\substack{ j=1\\
				\lambda \in \Lambda}}^{p}\left|\left\langle\mathbf{f}, \mathcal{D}_{2 N \lambda}\left(\mathbf{f}_{j}+\mathbf{g}_{j}\right)\right\rangle-\left\langle\mathbf{f},\mathcal{D}_{2 N \lambda} \mathbf{f}_{j}\right\rangle\right|^{2} } \notag\\
		& \geq \sqrt{\sum_{\substack{j=1 \\ \lambda \in \Lambda}}^{p} \left| \left\langle \mathbf{f}, \mathcal{D}_{2 N \lambda} \mathbf{f}_{j} \right\rangle \right|^{2}} - \sqrt{\sum_{\substack{j=1 \\ \lambda \in \Lambda}}^{p} \left| \left\langle \mathbf{f}, \mathcal{D}_{2 N \lambda} \left(\mathbf{f}_{j} + \mathbf{g}_{j} \right) \right\rangle \right|^{2}} \notag \\
		& \geq \sqrt{a_{0} \|\mathbf{f}\|^{2}} - \sqrt{\left(2^{p-1} \epsilon^{2} n^{2}\right) \|\mathbf{f}\|^{2}} \notag \\
		& = \left(\sqrt{a_{0}} - \sqrt{2^{p-1} \epsilon^{2} n^{2}} \right) \|\mathbf{f}\|, \  \mathbf{f} \in \ell^{2}\left(\Lambda,\mathcal{M}_{n}\right).
	\end{align}
	From inequalities \eqref{Per4.5} and \eqref{Per4.6}, we conclude that
$\left\{ \mathcal{D}_{2 N \lambda} \mathbf{g}_{j}\right\}_{\lambda\in \Lambda \atop j \in \mathbb{N}_p}$
	is a   matrix-valued DFNS for $\ell^{2}\left(\Lambda,\mathcal{M}_{n}\right)$ with the desired frame bounds.
	\endproof
	Now, we illustrate Theorem \ref{PerMI} with an example.
	\begin{example}
	Let $N=2$, $r=1$, $n=2 $ and $p=8$. Then, $\Lambda=\left\{0, \frac{1}{2}\right\}+2\mathbb{Z}$, $\Omega=\left[0, \frac{1}{2}\left)\cup\left[1, \frac{3}{2}\left),\right.\right.\right.\right.$ $j \in \mathbb{N}_{8}$. Consider the  matrix-valued DFNS $\left\{ \mathcal{D}_{2 N \lambda} \mathbf{f}_{j}\right\} _ {\lambda\in \Lambda \atop j \in [8]}$  as given in Example \ref{Ex3.10}. Let  $ a_{0} $ and $ b_{0}$ be its frame
For $j \in \mathcal{M}_{8}$, define $\mathbf{g}_{j}=\left\{\mathbf{g}_{j}(\lambda)\right\}_{\lambda \in \Lambda} \subset \ell^{2}\left(\Lambda,\mathcal{M}_{2}\right)$ as follows:

\begin{align*}
\begin{aligned}
	& \mathbf{g}_{1}(0)=\left[\begin{array}{cc}
		-\frac{24}{25}&0\\
		0&-\frac{24}{25}
	\end{array}\right], \quad \mathbf{g}_{1}(4)=\left[\begin{array}{cc}
		0&-1\\
		-1&0
	\end{array}\right], \quad \mathbf{g}_{1}(\lambda)=\left[\begin{array}{cc}
		0 & 0\\
		0&0
	\end{array}\right] \text { for } \lambda \in \Lambda \backslash\{0,4\}; \\
	& \mathbf{g}_{2}(0)=\left[\begin{array}{cc}
		-\frac{24}{25}&0 \\
		0&\frac{24}{25}
	\end{array}\right], \quad \mathbf{g}_{2}(4)=\left[\begin{array}{cc}
		0&i\\
		-i&0
	\end{array}\right], \quad \mathbf{g}_{2}(\lambda)=\left[\begin{array}{cc}
		0 &0\\
		0&0
	\end{array}\right] \text { for } \lambda \in \Lambda \backslash\{0,4\};\\
	&\mathbf{g}_{3}(0)=\left[\begin{array}{cc}
		0 & -\frac{24}{25}i\\
		-\frac{24}{25}i&0
	\end{array}\right], \quad \mathbf{g}_{3}(4)=\left[\begin{array}{cc}
		0 & -1\\
		-1& 0
	\end{array}\right], \quad \mathbf{g}_{3}(\lambda)=\left[\begin{array}{cc}
		0 & 0\\
		0& 0
	\end{array}\right] \text {for} \lambda \in \Lambda \backslash\{0,4\}; \\
	&\mathbf{g}_{4}(0)=\left[\begin{array}{cc}
		0& -\frac{24}{25}\\
		-\frac{24}{25}&0
	\end{array}\right], \quad \mathbf{g}_{4}(4)=\left[\begin{array}{cc}
		-1 & 0\\
		0&1
	\end{array}\right], \quad \mathbf{g}_{4}(\lambda)=\left[\begin{array}{cc}
		0 & 0\\
		0& 0
	\end{array}\right] \text { for } \lambda \in \Lambda \backslash\{0,4\}; \\
	& \mathbf{g}_{5}(\frac{1}{2})=\left[\begin{array}{cc}
		-1& 0\\
		0&-1
	\end{array}\right], \quad \mathbf{g}_{5}(\frac{1}{2}+4)=\left[\begin{array}{cc}
		0& -1\\
		-1 &0
	\end{array}\right], \quad \mathbf{g}_{5}(\lambda)=\left[\begin{array}{cc}
		0 & 0\\
		0& 0
	\end{array}\right] \text { for } \lambda \in \Lambda \backslash\{\frac{1}{2},\frac{1}{2}+4\}; \\
	& \mathbf{g}_{6}(\frac{1}{2})=\left[\begin{array}{cc}
		-1 & 0\\
		0&1
	\end{array}\right], \quad \mathbf{g}_{6}(\frac{1}{2}+4)=\left[\begin{array}{cc}
		0 & i\\
		-i &0
	\end{array}\right], \quad \mathbf{g}_{6}(\lambda)=\left[\begin{array}{cc}
		0 & 0\\
		0& 0
	\end{array}\right] \text { for } \lambda \in \Lambda \backslash\{\frac{1}{2},\frac{1}{2}+4\}; \\
	& \mathbf{g}_{7}(\frac{1}{2})=\left[\begin{array}{cc}
		0 & -i\\
		i&0
	\end{array}\right], \quad \mathbf{g}_{7}(\frac{1}{2}+4)=\left[\begin{array}{cc}
		0 & -1\\
		-1 &0
	\end{array}\right], \quad \mathbf{g}_{7}(\lambda)=\left[\begin{array}{cc}
		0 & 0\\
		0& 0
	\end{array}\right] \text { for } \lambda \in \Lambda \backslash\{\frac{1}{2},\frac{1}{2}+4\}; \\
	& \mathbf{g}_{8}(\frac{1}{2})=\left[\begin{array}{cc}
		0& -1\\
		-1&0
	\end{array}\right], \quad \mathbf{g}_{8}(\frac{1}{2}+4)=\left[\begin{array}{cc}
		-1 & 0\\
		0 &1
	\end{array}\right], \quad \mathbf{g}_{8}(\lambda)=\left[\begin{array}{cc}
		0 & 0\\
		0& 0
	\end{array}\right] \text { for } \lambda \in \Lambda \backslash\{\frac{1}{2},\frac{1}{2}+4\},
	\end{aligned}
\end{align*}
For $ j \in  \mathbb{N}_{8}$, the respective Fourier transforms  $\mathcal{F}\left(\mathbf{g}_{j}\right)(x)$ are given by
\begin{align*}
	&\mathcal{F}\left( \mathbf{g}_{1}\right)(x)=\left[\begin{array}{cc}
		-\frac{24}{25}& -e^{8\pi i x} \\
		-e^{8\pi i x}&-\frac{24}{25}
	\end{array}\right] ,&
	\mathcal{F}\left(\mathbf{g}_{2}\right)(x)=\left[\begin{array}{cc}
		-\frac{24}{25}& i e^{8\pi i x}\\
		-i e^{8\pi i x} &\frac{24}{25}
	\end{array}\right], \\
	&\mathcal{F}\left(\mathbf{g}_{3}\right)(x)=\left[\begin{array}{cc}
		0&-\frac{24}{25}i-e^{8\pi i x} \\
		-\frac{24}{25}i-e^{8 \pi ix} &0
	\end{array}\right],&
	\mathcal{F}\left(\mathbf{g}_{4}\right)(x)=\left[\begin{array}{cc}
		-e^{8\pi i x} & -\frac{24}{25}\\
		-\frac{24}{25}& e^{8 \pi i x}
	\end{array}\right],\\
	&\mathcal{F}\left( \mathbf{g}_{5}\right)(x)=\left[\begin{array}{cc}
		-e^{\pi i x} & -e^{9\pi i x} \\
		-e^{9\pi i x}&-e^{\pi i x}
	\end{array}\right] ,&
	\mathcal{F}\left(\mathbf{g}_{6}\right)(x)=\left[\begin{array}{cc}
		-e^{\pi i x}& ie^{9 \pi i x}\\
		-ie^{9\pi i x}  &e^{ \pi i x}
	\end{array}\right], \\
	&\mathcal{F}\left(\mathbf{g}_{7}\right)(x)=\left[\begin{array}{cc}
		0&-ie^{\pi i x} -e^{9\pi i x} \\
		ie^{ \pi ix}-e^{9\pi i x}  &0
	\end{array}\right],&
	\mathcal{F}\left(\mathbf{g}_{8}\right)(x)=\left[\begin{array}{cc}
		-e^{9 \pi i x} & -e^{\pi i x} \\
		-e^{\pi i x} & e^{9 \pi i x}
	\end{array}\right].
\end{align*}
For $ j \in \mathbb{N}_{8}$ and $ x \in \Omega $, it is straightforward to verify that
\begin{align*}
&\left\|\mathcal{F}\left(\mathbf{f}_{j} + \mathbf{g}_{j}\right)(x)\right\|_{\mathcal{M}_{2}} =
\left\| \mathcal{F}(\mathbf{f}_{j}) (x) + \mathcal{F}\left(\mathbf{g}_{j}\right)(x) \right\| \leq \frac{1}{25} =\epsilon \, (\text{say}),
\text{ and }\\
&\epsilon<2^{p-1} \epsilon^{2} n^{2} = 2^{7} \times \epsilon^{2} \times 4 =  0.8192 < 1= a_{0}.
\end{align*}
Thus, by  Theorem \ref{PerMI}, the collection $\left\{ \mathcal{D}_{2 N \lambda} \mathbf{g}_{j}\right\}_{\lambda\in \Lambda \atop j \in \mathbb{N}_{8}}$
is a matrix-valued DFNS for $ \ell^{2}\left(\Lambda,\mathcal{M}_{2}\right)$.
\end{example}
We recall that errors are generally measured using relative error \cite{Ford}.  Relative error quantifies how far an estimated value is from the true value, relative to the magnitude of the true value. Mathematically, it is defined as follows:
\begin{align*}
\text{Relative Error} = \frac{|\hat{\theta} - \theta|}{|\theta|},
\end{align*}
where \( \hat{\theta} \) is the estimated value and \( \theta \) is the true (reference) value. This measure is especially useful when comparing errors across values of different scales.

In statistics or analysis, relative quantities can be {less sensitive to outliers or noise.\\
For example, in the expression:
\[
\frac{\left\| \mathcal{F}(\mathbf{g}_j)(x) - \mathcal{F}(\mathbf{f}_j)(x) \right\|}{\left\| \mathcal{F}(\mathbf{f}_j)(x) \right\|}, \, j \in \mathbb{N}_{p},
\]
we are not just measuring the magnitude of the perturbation, but rather the perturbation  \emph{relative to the size of the signal}. This often gives a more meaningful and stable measure of deviation, especially when comparing across different scales or conditions. Motivated by the notion of relative error, the following theorem gives frame conditions for perturbed sequences.

\begin{theorem}\label{PerMIIx}
Suppose
\begin{enumerate}
  \item $\left\{ \mathcal{D}_{2 N \lambda} \mathbf{f}_{j}\right\}_{\lambda\in \Lambda \atop j \in \mathbb{N}_{p}}$  is a  matrix-valued DFNS for $\ell^{2}\left(\Lambda, \mathcal{M}_n\right)$ with frames bounds $a_{0}$ and $b_{0}$.
  \item  $\left\{\mathbf{g}_{j}\right\}_{j \in \mathbb{N}_{p}} \subset \ell^{2}\left(\Lambda,\mathcal{M}_{n}\right)$ satisfies
$\frac{\left\| \mathcal{F}(\mathbf{g}_{j})(x)- \mathcal{F}(\mathbf{f}_{j})(x) \right\|_{\mathcal{M}_n}}{\left\| \mathcal{F}(\mathbf{f}_{j})(x) \right\|_{\mathcal{M}_n}} \leq \epsilon$, where $\epsilon > 0$ with $\epsilon^2 < \frac{a_0}{2^{p-1} (N+b_0)^{2} n^2}$ and $x \in \Omega$.\label{per2}
\end{enumerate}
Then, the collection
$\left\{ \mathcal{D}_{2 N \lambda} \mathbf{g}_{j} \right\}_{\lambda\in \Lambda \atop j \in \mathbb{N}_{8}}$ is also a  matrix-valued DFNS for
$\ell^{2}\left(\Lambda ,\mathcal{M}_n\right)$ with frame bounds
$\left( \sqrt{a_{0}} - \sqrt{2^{p-1} \epsilon^{2} (N+ b_{0})^{2} n^{2}}\right)^{2}$ and $ \left(2^{p}\epsilon^{2}(N+ b_{0})^{2}n^{2}+2b_{0}\right)$.
\end{theorem}	
\proof
Using  Lemma \ref{Blem1}, for each  $\mathbf{f} \in \ell^{2}\left(\Lambda ,\mathcal{M}_n\right)$, we have
		 \begin{align*}
		\sum_{\substack{ j=1\\
		 			\lambda \in \Lambda}}^{p}\left|\left\langle\mathbf{f}, \mathcal{D}_{2 N \lambda} \mathbf{g}_{j}\right\rangle\right|^{2}
		 	= & \sum_{\substack{ j=1\\
		 			\lambda \in \Lambda}}^{p}\left|\left\langle\mathbf{f}, \mathcal{D}_{2 N \lambda}\left(\mathbf{g}_{j}-\mathbf{f}_{j}\right)\right\rangle+\left\langle\mathbf{f},\mathcal{D}_{2 N \lambda} \mathbf{f}_{j}\right\rangle\right|^{2}  \\
		 	\leq & 2\left(\sum_{\substack{ j=1 \\
		 			\lambda \in \Lambda}}^{p}\left|\left\langle\mathbf{f}, \mathcal{D}_{2 N \lambda}\left(\mathbf{g}_{j}-\mathbf{f}_{j}\right)\right\rangle\right|^{2}+\sum_{\substack{j=1 \\ \lambda \in \Lambda}}^{p}\left|\left\langle\mathbf{f},\mathcal{D}_{2 N \lambda} \mathbf{f}_{j}\right\rangle\right|^{2}\right).
		 \end{align*}
Using hypothesis \eqref{per2}, and corollary \ref{corM3} , we get
\begin{align*}
	\left\| \mathcal{F}(\mathbf{g}_{j} - \mathbf{f}_{j})(x) \right\|_{\mathcal{M}_n} \leq \epsilon \left\| \mathcal{F}(\mathbf{f}_{j}(x)) \right\|_{\mathcal{M}_n}\leq \epsilon (N+b_{0}).
\end{align*}
Now,  Theorem \ref{thmM1} gives
	\begin{equation*}
	\sum_{\substack{j=1 \\  \lambda \in \Lambda}}^{p}\left|\left\langle\mathbf{f},  \mathcal{D}_{2 N \lambda}\left(\mathbf{g}_{j}-\mathbf{f}_{j}\right)\right\rangle\right|^{2} \leq 2^{p-1} \epsilon^{2}(N+b_{0})^{2}n^{2}\|\mathbf{f}\|^{2},  \, \mathbf{f} \in \ell^{2}\left(\Lambda ,\mathcal{M}_n\right).
\end{equation*}
Therefore, we  have obtain  the upper frame inequality
		 \begin{equation}\label{eqii1}
		 	\sum_{\substack{j=1 \\ \lambda \in \Lambda}}^{p}\left|\left\langle\mathbf{f}, \mathcal{D}_{2 N \lambda} \mathbf{g}_{j}\right\rangle\right|^{2} \leq\left(2^{p}\epsilon^{2}(N+ b_{0})^{2}n^{2}+2b_{0}\right)\|\mathbf{f}\|^{2}, \, \mathbf{f} \in \ell^{2}\left(\Lambda ,\mathcal{M}_n\right).
		 \end{equation}
For the lower frame inequality,  we compute
		 \begin{align}\label{eqaa1}
		 	 \sqrt{\sum_{\substack{j=1 \\ \lambda \in \Lambda}}^{p} \left| \left\langle \mathbf{f}, \mathcal{D}_{2 N \lambda} \mathbf{g}_{j} \right\rangle \right|^{2}}
		 	= & \sqrt{\sum_{\substack{ j=1\\
		 	\lambda \in \Lambda}}^{p}\left|\left\langle\mathbf{f}, \mathcal{D}_{2 N \lambda}\left(\mathbf{g}_{j}-\mathbf{f}_{j}\right)\right\rangle+\left\langle\mathbf{f},\mathcal{D}_{2 N \lambda} \mathbf{f}_{j}\right\rangle\right|^{2} }\notag \\
		 	& \geq \sqrt{\sum_{\substack{j=1 \\ \lambda \in \Lambda}}^{p} \left| \left\langle \mathbf{f}, \mathcal{D}_{2 N \lambda} \mathbf{f}_{j} \right\rangle \right|^{2}} - \sqrt{\sum_{\substack{j=1 \\ \lambda \in \Lambda}}^{p} \left| \left\langle \mathbf{f}, \mathcal{D}_{2 N \lambda} \left(\mathbf{g}_{j} -\mathbf{f}_{j} \right) \right\rangle \right|^{2}} \notag \\
		& \geq \sqrt{a_{0} \|\mathbf{f}\|^{2}} - \sqrt{2^{p-1} \epsilon^{2} (N+ b_{0})^{2}n^{2} \|\mathbf{f}\|^{2}} \notag \\
		 	& = \left( \sqrt{a_{0}} - \sqrt{2^{p-1} \epsilon^{2} (N +b_{0})^{2} n^{2}} \right) \|\mathbf{f}\| \, \text{for all} \  \mathbf{f} \in \ell^{2}\left(\Lambda, \mathcal{M}_{n}\right).
		 \end{align}
Hence, from \eqref{eqii1} and \eqref{eqaa1}, we conclude that the collection $\left\{ \mathcal{D}_{2 N \lambda} \mathbf{g}_{j} \right\}_{\lambda \in \Lambda \atop j \in \mathbb{N}_{p}}$ form a	matrix-valued DFNS for $\ell^{2}\left(\Lambda ,\mathcal{M}_n\right)$ with the desired frame bounds.
		\endproof

$$\text{\textbf{Funding}}$$
 The  research of Manisha Chhillar  is  supported by the  Council of Scientific and Industrial Research (CSIR) under the file number 09/0045(16423)/2023-EMR-I.

\textbf{Hari Krishan Malhotra}, Department of Mathematics,
	 Manav Rachna Internatioanl Institute of Research and Studies, Faridabad, Haryana, India\\
\emph{Email:} maths.hari67@gmail.com\\

\textbf{Manisha Chhillar}, Department of Mathematics,
Department of Mathematics,		
University of Delhi, Delhi-110007.	\\
\emph{Email:} manshachhilar71@gmail.com\\

\textbf{Lalit Kumar Vashisht}, Department of Mathematics,
		University of Delhi, Delhi-110007.\\
Email: lalitkvashisht@gmail.com

\end{document}